\documentclass[journal]{IEEEtran}
\pdfoutput=1

\usepackage{cite}
\usepackage{graphicx}
\usepackage[breaklinks=false, allcolors=black, colorlinks=true]{hyperref}

\usepackage{amsfonts}
\usepackage{amsmath} 
\usepackage{amssymb}
\usepackage{amsthm}	
\usepackage{enumerate}
\usepackage{mathrsfs}
\usepackage{booktabs}
\usepackage{epstopdf}
\usepackage[percent]{overpic}
\usepackage{bm}

\makeatletter

\newcommand{\Rmnum}[1]{\expandafter\@slowromancap\romannumeral #1@}
\makeatother

\theoremstyle{plain} 
\newtheorem{theore}{Theorem}
\newtheorem{corollar}{Corollary}
\newtheorem{lemm}[theore]{Lemma}
\newtheorem{definitio}{Definition}
\newtheorem{proble}{Problem}
\newtheorem{assumptio}{Assumption}
\newtheorem{propert}{Property}
\newtheorem{propositio}{Proposition}
\newtheorem{conjectur}{Conjecture}
\newtheorem{exampl}{Example}
\newtheorem{rul}{Rule}

\newtheorem{remar}{Remark}
\newtheorem{fac}{Fact}

\def\nn{\nonumber}
\def\beq{\begin{equation}}
\def\eeq{\end{equation}}
\def\beqs{\begin{equation*}}
\def\eeqs{\end{equation*}}
\def\balign{\begin{align}}
\def\ealign{\end{align}}
\def\bea{\begin{eqnarray}}
\def\eea{\end{eqnarray}}
\def\ba{\begin{array}}
\def\ea{\end{array}}
\def\bcenter{\begin{center}}
\def\ecenter{\end{center}}
\def\bitem{\begin{itemize}}
\def\eitem{\end{itemize}}
\def\benum{\begin{enumerate}}
\def\eenum{\end{enumerate}}
\def\bdoc{\begin{document}}
\def\edoc{\end{document}}
\def\defeq{{\stackrel{\Delta}{=}}}
\def\nn{\nonumber}

\def\eg{{e.g., \/}}
\def\etal{{\it et al. \/}}
\def\viz{{\iet viz.,\ \/}}
\def\ie{{i.e.,\ \/}}
\def\etc{{\it etc.,\ \/}}
\def\cf{{cf.,\ \/}}

\usepackage[usenames,dvipsnames,svgnames,table]{xcolor}

\def\black{\textcolor{black}}
\def\yellow{\textcolor{yellow}}
\def\green{\textcolor{green}}
\def\tcbg{\textcolor{bgrd}}
\def\magenta{\textcolor{magenta}}
\def\white{\textcolor{white}}
\def\gray{\textcolor{gray}}

\newcommand\Quote[1]{\lq\textsl{#1}\rq}
\newcommand\fr[2]{{\textstyle\frac{#1}{#2}}}

\newcommand{\Ambb}{\mathbb{A}}
\newcommand{\Bmbb}{\mathbb{B}}
\newcommand{\Cmbb}{\mathbb{C}}
\newcommand{\Dmbb}{\mathbb{D}}
\newcommand{\Embb}{\mathbb{E}}
\newcommand{\Fmbb}{\mathbb{F}}
\newcommand{\Gmbb}{\mathbb{G}}
\newcommand{\Hmbb}{\mathbb{H}}
\newcommand{\Imbb}{\mathbb{I}}
\newcommand{\Jmbb}{\mathbb{J}}
\newcommand{\Kmbb}{\mathbb{K}}
\newcommand{\Lmbb}{\mathbb{L}}
\newcommand{\Mmbb}{\mathbb{M}}
\newcommand{\Nmbb}{\mathbb{N}}
\newcommand{\Ombb}{\mathbb{O}}
\newcommand{\Pmbb}{\mathbb{P}}
\newcommand{\Qmbb}{\mathbb{Q}}
\newcommand{\Rmbb}{\mathbb{R}}
\newcommand{\Smbb}{\mathbb{S}}
\newcommand{\Tmbb}{\mathbb{T}}
\newcommand{\Umbb}{\mathbb{U}}
\newcommand{\Vmbb}{\mathbb{V}}
\newcommand{\Wmbb}{\mathbb{W}}
\newcommand{\Xmbb}{\mathbb{X}}
\newcommand{\Ymbb}{\mathbb{Y}}
\newcommand{\Zmbb}{\mathbb{Z}}

\newcommand{\Amsc}{\mathcal{A}}
\newcommand{\Bmsc}{\mathcal{B}}
\newcommand{\Cmsc}{\mathcal{C}}
\newcommand{\Dmsc}{\mathcal{D}}
\newcommand{\Emsc}{\mathcal{E}}
\newcommand{\Fmsc}{\mathcal{F}}
\newcommand{\Gmsc}{\mathcal{G}}
\newcommand{\Hmsc}{\mathcal{H}}
\newcommand{\Imsc}{\mathcal{I}}
\newcommand{\Jmsc}{\mathcal{J}}
\newcommand{\Kmsc}{\mathcal{K}}
\newcommand{\Lmsc}{\mathcal{L}}
\newcommand{\Mmsc}{\mathcal{M}}
\newcommand{\Nmsc}{\mathcal{N}}
\newcommand{\Omsc}{\mathcal{O}}
\newcommand{\Pmsc}{\mathcal{P}}
\newcommand{\Qmsc}{\mathcal{Q}}
\newcommand{\Rmsc}{\mathcal{R}}
\newcommand{\Smsc}{\mathcal{S}}
\newcommand{\Tmsc}{\mathcal{T}}
\newcommand{\Umsc}{\mathcal{U}}
\newcommand{\Vmsc}{\mathcal{V}}
\newcommand{\Wmsc}{\mathcal{W}}
\newcommand{\Xmsc}{\mathcal{X}}
\newcommand{\Ymsc}{\mathcal{Y}}
\newcommand{\Zmsc}{\mathcal{Z}}

\newcommand{\Ascr}{\mathscr{A}}
\newcommand{\Bscr}{\mathscr{B}}
\newcommand{\Cscr}{\mathscr{C}}
\newcommand{\Dscr}{\mathscr{D}}
\newcommand{\Escr}{\mathscr{E}}
\newcommand{\Fscr}{\mathscr{F}}
\newcommand{\Gscr}{\mathscr{G}}
\newcommand{\Hscr}{\mathscr{H}}
\newcommand{\Iscr}{\mathscr{I}}
\newcommand{\Jscr}{\mathscr{J}}
\newcommand{\Kscr}{\mathscr{K}}
\newcommand{\Lscr}{\mathscr{L}}
\newcommand{\Mscr}{\mathscr{M}}
\newcommand{\Nscr}{\mathscr{N}}
\newcommand{\Oscr}{\mathscr{O}}
\newcommand{\Pscr}{\mathscr{P}}
\newcommand{\Qscr}{\mathscr{Q}}
\newcommand{\Rscr}{\mathscr{R}}
\newcommand{\Sscr}{\mathscr{S}}
\newcommand{\Tscr}{\mathscr{T}}
\newcommand{\Uscr}{\mathscr{U}}
\newcommand{\Vscr}{\mathscr{V}}
\newcommand{\Wscr}{\mathscr{W}}
\newcommand{\Xscr}{\mathscr{X}}
\newcommand{\Yscr}{\mathscr{Y}}
\newcommand{\Zscr}{\mathscr{Z}}

\def\Atiny{\mbox{\tiny{A}}}
\def\Btiny{\mbox{\tiny{B}}}
\def\Ctiny{\mbox{\tiny{C}}}
\def\Dtiny{\mbox{\tiny{D}}}
\def\Etiny{\mbox{\tiny{E}}}
\def\Ctiny{\mbox{\tiny{C}}}
\def\Ftiny{\mbox{\tiny{F}}}
\def\Gtiny{\mbox{\tiny{G}}}
\def\Htiny{\mbox{\tiny{H}}}
\def\Itiny{\mbox{\tiny{I}}}
\def\Jtiny{\mbox{\tiny{J}}}
\def\Ktiny{\mbox{\tiny{K}}}
\def\Ltiny{\mbox{\tiny{L}}}
\def\Mtiny{\mbox{\tiny{M}}}
\def\Ntiny{\mbox{\tiny{N}}}
\def\Otiny{\mbox{\tiny{O}}}
\def\Ptiny{\mbox{\tiny{P}}}
\def\Qtiny{\mbox{\tiny{Q}}}
\def\Rtiny{\mbox{\tiny{R}}}
\def\Stiny{\mbox{\tiny{S}}}
\def\Ttiny{\mbox{\tiny{T}}}
\def\Utiny{\mbox{\tiny{U}}}
\def\Vtiny{\mbox{\tiny{V}}}
\def\Wtiny{\mbox{\tiny{W}}}
\def\Xtiny{\mbox{\tiny{X}}}
\def\Ytiny{\mbox{\tiny{Y}}}
\def\Ztiny{\mbox{\tiny{Z}}}

\def\alphabf{\hbox{\boldmath$\alpha$\unboldmath}}
\def\betabf{\hbox{\boldmath$\beta$\unboldmath}}
\def\gammabf{\hbox{\boldmath$\gamma$\unboldmath}}
\def\deltabf{\hbox{\boldmath$\delta$\unboldmath}}
\def\epsilonbf{\hbox{\boldmath$\epsilon$\unboldmath}}
\def\zetabf{\hbox{\boldmath$\zeta$\unboldmath}}
\def\etabf{\hbox{\boldmath$\eta$\unboldmath}}
\def\iotabf{\hbox{\boldmath$\iota$\unboldmath}}
\def\kappabf{\hbox{\boldmath$\kappa$\unboldmath}}
\def\lambdabf{\hbox{\boldmath$\lambda$\unboldmath}}
\def\mubf{\hbox{\boldmath$\mu$\unboldmath}}
\def\nubf{\hbox{\boldmath$\nu$\unboldmath}}
\def\xibf{\hbox{\boldmath$\xi$\unboldmath}}
\def\pibf{\hbox{\boldmath$\pi$\unboldmath}}
\def\rhobf{\hbox{\boldmath$\rho$\unboldmath}}
\def\sigmabf{\hbox{\boldmath$\sigma$\unboldmath}}
\def\taubf{\hbox{\boldmath$\tau$\unboldmath}}
\def\upsilonbf{\hbox{\boldmath$\upsilon$\unboldmath}}
\def\phibf{\hbox{\boldmath$\phi$\unboldmath}}
\def\chibf{\hbox{\boldmath$\chi$\unboldmath}}
\def\psibf{\hbox{\boldmath$\psi$\unboldmath}}
\def\omegabf{\hbox{\boldmath$\omega$\unboldmath}}
\def\Sigmabf{\hbox{$\bf \Sigma$}}
\def\Upsilonbf{\hbox{$\bf \Upsilon$}}
\def\Omegabf{\hbox{$\bf \Omega$}}
\def\Deltabf{\hbox{$\bf \Delta$}}
\def\Gammabf{\hbox{$\bf \Gamma$}}
\def\Thetabf{\hbox{$\bf \Theta$}}
\def\Lambdabf{\mbox{$ \bf \Lambda $}}
\def\Xibf{\hbox{\bf$\Xi$}}
\def\Pibf{{\bf \Pi}}
\def\thetabf{{\mbox{\boldmath$\theta$\unboldmath}}}
\def\Upsilonbf{\hbox{\boldmath$\Upsilon$\unboldmath}}
\newcommand{\Phibf}{\mbox{${\bf \Phi}$}}
\newcommand{\Psibf}{\mbox{${\bf \Psi}$}}

\def\abf{{\bf a}}
\def\bbf{{\bf b}}
\def\cbf{{\bf c}}
\def\dbf{{\bf d}}
\def\ebf{{\bf e}}
\def\fbf{{\bf f}}
\def\gbf{{\bf g}}
\def\hbf{{\bf h}}
\def\ibf{{\bf i}}
\def\jbf{{\bf j}}
\def\kbf{{\bf k}}
\def\lbf{{\bf l}}
\def\mbf{{\bf m}}
\def\nbf{{\bf n}}
\def\obf{{\bf o}}
\def\pbf{{\bf p}}
\def\qbf{{\bf q}}
\def\rbf{{\bf r}}
\def\sbf{{\bf s}}
\def\tbf{{\bf t}}
\def\ubf{{\bf u}}
\def\vbf{{\bf v}}
\def\wbf{{\bf w}}
\def\xbf{{\bf x}}
\def\ybf{{\bf y}}
\def\zbf{{\bf z}}
\def\rbf{{\bf r}}
\def\xbf{{\bf x}}
\def\ybf{{\bf y}}
\def\Abf{{\bf A}}
\def\Bbf{{\bf B}}
\def\Cbf{{\bf C}}
\def\Dbf{{\bf D}}
\def\Ebf{{\bf E}}
\def\Fbf{{\bf F}}
\def\Gbf{{\bf G}}
\def\Hbf{{\bf H}}
\def\Ibf{{\bf I}}
\def\Jbf{{\bf J}}
\def\Kbf{{\bf K}}
\def\Lbf{{\bf L}}
\def\Mbf{{\bf M}}
\def\Nbf{{\bf N}}
\def\Obf{{\bf O}}
\def\Pbf{{\bf P}}
\def\Qbf{{\bf Q}}
\def\Rbf{{\bf R}}
\def\Sbf{{\bf S}}
\def\Tbf{{\bf T}}
\def\Ubf{{\bf U}}
\def\Vbf{{\bf V}}
\def\Wbf{{\bf W}}
\def\Xbf{{\bf X}}
\def\Ybf{{\bf Y}}
\def\Zbf{{\bf Z}}
\def\Ac{{\cal A}}
\def\Bc{{\cal B}}
\def\Cc{{\cal C}}
\def\Dc{{\cal D}}
\def\Ec{{\cal E}}
\def\Fc{{\cal F}}
\def\Gc{{\cal G}}
\def\Hc{{\cal H}}
\def\Ic{{\cal I}}
\def\Jc{{\cal J}}
\def\Kc{{\cal K}}
\def\Lc{{\cal L}}
\def\Mc{{\cal M}}
\def\Nc{{\cal N}}
\def\Oc{{\cal O}}
\def\Pc{{\cal P}}
\def\Qc{{\cal Q}}
\def\Rc{{\cal R}}
\def\Sc{{\cal S}}
\def\Tc{{\cal T}}
\def\Uc{{\cal U}}
\def\Vc{{\cal V}}
\def\Wc{{\cal W}}
\def\Xc{{\cal X}}
\def\Yc{{\cal Y}}
\def\Zc{{\cal Z}}

\def\0bf{\mathbf{0}}
\def\1bf{\mathbf{1}}

\def\lambdabf{\boldsymbol{\lambda}}
\def\mubf{\boldsymbol{\mu}}
\def\gammabf{\boldsymbol{\gamma}}
\def\xibf{\boldsymbol{\xi}}
\def\sigmabf{\boldsymbol{\sigma}}
\def\zetabf{\boldsymbol{\zeta}}

\def\Lambdabf{\boldsymbol{\Lambda}}
\def\Mubf{\boldsymbol{\Mu}}
\def\Gammabf{\boldsymbol{\Gamma}}
\def\Xibf{\boldsymbol{\Xi}}
\def\Sigmabf{\boldsymbol{\Sigma}}
\def\Zetabf{\boldsymbol{\Zeta}}

\newcommand\independent{\protect\mathpalette{\protect\independenT}{\perp}}
\def\independenT#1#2{\mathrel{\rlap{$#1#2$}\mkern2mu{#1#2}}}

\newcommand{\bsf}[1]{\textsf{\textbf{#1}}}
\newcommand{\lbsf}[1]{\textsf{\large  \textbf{#1}}}
\newcommand{\Lbsf}[1]{\textsf{\Large  \textbf{#1}}}
\newcommand{\hbsf}[1]{\textsf{\huge  \textbf{#1}}}

\newcommand{\myminipage}[3]{\begin{minipage}[#1]{#2}{#3} \end{minipage}}
\newcommand{\sbs}[4]{\myminipage{c}{#1}{#3} \hfill
\myminipage{c}{#2}{#4}}

\newcommand{\myfig}[2]{\centerline{\psfig{figure=#1,width=#2,silent=}}}
\newcommand{\myfigh}[2]{\centerline{\psfig{figure=#1,height=#2,silent=}}}
\newcommand{\myfigwh}[3]{\centerline{\psfig{figure=#1,width=#2,height=#3,silent=}}}

\newcommand{\beqa}{\begin{eqnarray}}
\newcommand{\eeqa}{\end{eqnarray}}
\newcommand{\beqan}{\begin{eqnarray*}}
\newcommand{\eeqan}{\end{eqnarray*}}
\newcommand{\dst}[1]{\displaystyle{ #1 }}
\newcounter{pf}

\newcommand{\smax}[1] { \bar \sigma \left( #1 \right) }
\newcommand{\Rn}{{\mathbb R}^n}
\newcommand{\R}{{\mathbb R}}
\newcommand{\C}{{\mathbb C}}
\newcommand{\Rm}{\mathbb{R}^m}
\newcommand{\Rmn}{\mathbb{R}^{m \times n}}
\newcommand{\Rpq}{\mathbb{R}^{p \times q}}
\newcommand{\Cn}{\mathbb{C}^n}
\newcommand{\Cm}{\mathbb{C}^m}
\newcommand{\Cnn}{\mathbb{C}^{n \times n}}
\newcommand{\Cmn}{\mathbb{C}^{m \times n}}
\newcommand{\ip}[1]{\left\langle #1 \right\rangle}
\newcommand{\rank}{\mbox{rank}}
\newcommand{\Span}{\mbox{\rm Span }}
\newcommand{\Trace}{\mbox{\rm Tr }}
\newcommand{\Spec}{\mbox{\rm Spec }}
\newcommand{\vectornorm}[1]{\left\|#1\right\|}

\newcommand{\pd}[2]{\frac{\partial #1}{\partial #2}}
\newcommand{\ppd}[3]{\frac{\partial^2 #1}{\partial #2 \partial #3}}

\newcommand{\thtilde}{\tilde{\theta}}
\newcommand{\thnom}{\theta^\circ}
\newcommand{\thopt}{\theta^{\mbox{\small opt}}}
\newcommand{\thhat}{{\hat{\theta}}}
\newcommand{\Tho}{\Theta^\circ}
\newcommand{\tho}{\theta^\circ}
\newcommand{\np}{{n_p}}

\newcommand{\ii}{{[i]}}
\newcommand{\II}{{[i+1]}}
\newcommand{\iii}{{[ii]}}
\newcommand{\jj}{{[j]}}
\newcommand{\kk}{{[k]}}
\newcommand{\thi}{{\theta^\ii}}
\newcommand{\thI}{{\theta^\II}}
\newcommand{\di}{{d^\ii}}
\newcommand{\gi}{{g^\ii}}
\newcommand{\Hi}{{\HH^\ii}}
\newcommand{\thK}{\theta^{(k+1)}}
\newcommand{\gk}{{g^{(k)}}}
\newcommand{\Hk}{{{\cal H}^{(k)}}}

\newcommand{\bfdelta}{{\bf \Delta}}

\newcommand{\Exp}[1]{\exp \left\{ #1 \right\}} 
\newcommand{\gaussian}[1]{\mathbb{N} \left( #1 \right)}
\newcommand{\uniform}[1]{\mathbb{U} \left[ #1 \right]}
\newcommand{\exponential}[1]{\mathbb{E} \left[ #1 \right]}
\newcommand{\EXP}[1]{\EEXP \left[ #1 \right]} 
\newcommand{\EEXP}{\mbox{\bsf{E}}} 
\newcommand{\Prob}[1]{\mbox{{\sf Pr}} \left(#1 \right)}
\newcommand{\convas}{\stackrel{as}{\longrightarrow}}
\newcommand{\convinp}{\stackrel{p}{\longrightarrow}}
\newcommand{\convind}{\stackrel{d}{\longrightarrow}}
\newcommand{\convqm}{\stackrel{qm}{\longrightarrow}}
\newcommand{\sss}[1]{{_{#1}}}
\newcommand{\density}[2]{p_{_{_{#1}}}\!\!\left(#2 \right)} 
\newcommand{\distro}[2]{P_{_{_{#1}}}\!\!\left(#2 \right)} 
\newcommand{\rxx}[1]{R_{_{#1}}\!} 
\newcommand{\sxx}[1]{S_{_{#1}}} 
\newcommand{\cov}[1]{\Lambda_{_{#1}}} 
\newcommand{\mean}[1]{m_{_{#1}}} 
\newcommand{\LS}[1]{\hat{#1}_{_{LS}}} 
\newcommand{\MV}[1]{\hat{#1}_{_{MV}}} 
\newcommand{\LMV}[1]{\hat{#1}_{_{LMV}}} 
\newcommand{\ML}[1]{\hat{#1}_{_{ML}}} 

\renewcommand{\arraystretch}{0.9}
\newcommand{\bmat}[1]{ \begin{bmatrix} #1 \end{bmatrix}}
\newcommand{\mat}[1]{ \left[ \begin{array}{cccccccc} #1 \end{array}
\right] }
\newcommand{\smallmat}[1]{\small{\mat{#1}}}
\newcommand{\sysblk}[4]{\begin{array}{c|cccc}#1&#2\\ \hline#3&#4
\end{array}}
\newcommand{\sysmat}[4]{\left[\sysblk{#1}{#2}{#3}{#4}\right]}
\newcommand{\SGeq}{\succ}
\newcommand{\SLeq}{\prec}
\newcommand{\Geq}{\succeq}
\newcommand{\Leq}{\preceq}

\newcommand{\Bset}{\mathbb{B}}
\newcommand{\Cset}{\mathbb{C}}
\newcommand{\Fset}{\mathbb{F}}
\newcommand{\Mset}{\mathbb{M}}
\newcommand{\Nset}{\mathbb{N}}
\newcommand{\Qset}{\mathbb{Q}}
\newcommand{\Rset}{\mathbb{R}}
\newcommand{\Sset}{\mathbb{S}}
\newcommand{\Tset}{\mathbb{T}}
\newcommand{\Uset}{\mathbb{U}}
\newcommand{\Vset}{\mathbb{V}}
\newcommand{\Wset}{\mathbb{W}}
\newcommand{\Zset}{\mathbb{Z}}

\newcommand{\Ical}{{\cal I}}
\newcommand{\Acal}{{\cal A}}
\newcommand{\Bcal}{{\cal B}}
\newcommand{\Ccal}{{\cal C}}
\newcommand{\Dcal}{{\cal D}}
\newcommand{\Ecal}{{\cal E}}
\newcommand{\Fcal}{{\cal F}}
\newcommand{\Gcal}{{\cal G}}
\newcommand{\Hcal}{{\cal H}}
\newcommand{\Kcal}{{\cal K}}
\newcommand{\Lcal}{{\cal L}}
\newcommand{\Mcal}{{\cal M}}
\newcommand{\Ncal}{{\cal N}}
\newcommand{\Pcal}{{\cal P}}
\newcommand{\Qcal}{{\cal Q}}
\newcommand{\Rcal}{{\cal R}}
\newcommand{\Scal}{{\cal S}}
\newcommand{\Tcal}{{\cal T}}
\newcommand{\Wcal}{{\cal W}}
\newcommand{\Ucal}{{\cal U}}
\newcommand{\Vcal}{{\cal V}}
\newcommand{\Xcal}{{\cal X}}
\newcommand{\Zcal}{{\cal Z}}

\newcommand{\FF}{{\bf F}}
\newcommand{\GG}{{\bf G}}
\newcommand{\HH}{{\bf H}}
\newcommand{\LL}{{\bf L}}
\newcommand{\NN}{{\bf N}}
\newcommand{\MM}{{\bf M}}
\newcommand{\PP}{{\bf P}}
\newcommand{\QQ}{{\bf Q}}
\newcommand{\RR}{{\bf R}}
\renewcommand{\SS}{{\bf S}}
\newcommand{\TT}{{\bf T}}
\newcommand{\VV}{{\bf V}}
\newcommand{\WW}{{\bf W}}

\newcommand{\thk}{\theta^{(k)}}
\newcommand{\thb}{\theta^{\rm opt}}
\newcommand{\alb}{\alpha^{\rm opt}}
\newcommand{\dk}{d^{(k)}}
\newcommand{\Hinf}{{\cal H}_\infty}
\newcommand{\Htwo}{{\cal H}_2}

\renewcommand{\arraystretch}{1.1}

\newcommand{\red}[1]{{\color{red} #1}}
\newcommand{\blue}[1]{#1}
\newcommand{\blueN}[1]{#1}
\newcommand{\blueF}[1]{#1}

\newcounter{l1}
\newcounter{l2}
\newcounter{l3}
\newcommand{\bdotlist}{\begin{list}{$\bullet$}{}}
\newcommand{\bboxlist}{\begin{list}{$\Box$}{}}
\newcommand{\bbboxlist}{\begin{list}{\raisebox{.005in}{{\tiny
$\blacksquare$ \ \ }}}{}}
\newcommand{\bdashlist}{\begin{list}{$-$}{} }
\newcommand{\blist}{\begin{list}{}{} }
\newcommand{\barablist}{\begin{list}{\arabic{l1}}{\usecounter{l1}}}
\newcommand{\balphlist}{\begin{list}{(\alph{l2})}{\usecounter{l2}}}
\newcommand{\bAlphlist}{\begin{list}{\Alph{l2}.}{\usecounter{l2}}}
\newcommand{\bdiamlist}{\begin{list}{$\diamond$}{}}
\newcommand{\bromalist}{\begin{list}{(\roman{l3})}{\usecounter{l3}}}

\newcommand{\thm}[1]{\noindent \begin{theorem} #1   \end{theorem}}
\newcommand{\prop}[1]{\begin{proposition} #1 \end{proposition}}
\newcommand{\lem}[1]{\begin{lemma} #1  \hfill $\blacksquare$ \end{lemma}}
\newcommand{\ex}[1]{\begin{example} {\rm #1} \end{example}}
\newcommand{\prf}[1]{ \noindent {\em Proof:} \, #1 \hfill $\blacksquare$}
\newcommand{\rem}[1]{\begin{remark} {\rm #1} \hfill $\Box$ \end{remark}}
\newcommand{\defn}[1]{\begin{definition} {\rm #1 } \end{definition}}
\newcommand{\prob}[1]{\begin{exercise} {\rm  #1 } \end{exercise}}
\newcommand{\cor}[1]{\begin{corollary}   #1  \end{corollary}}

\newcommand{\argmin}{\mathop{\rm argmin}}
\newcommand{\argmax}{\mathop{\rm argmax}}
\newcommand{\diag}{\mathop{\mathrm{diag}}}
\newcommand{\tr}{\mathop{\rm Tr}}
\newcommand{\conv}{\mathop{\rm conv}}
\newcommand{\var}{\mathop{\rm Var}}
\renewcommand{\b}[1]{\ensuremath{\boldsymbol{\mathrm{#1}}}}
\newcommand{\ms}{{\rm MS}}
\newcommand{\tcs}{{\rm TCS}}
\newcommand{\scs}{{\rm SCS}}
\newcommand{\E}[1]{\b{\mu}_{{#1}}}
\newcommand{\Var}[1]{{\Sigma_{#1}}}
\newcommand{\newbus}{fourteen}

\IEEEoverridecommandlockouts

\renewcommand{\Exp}{\mathbb{E}}
\renewcommand{\Cn}{\mathbf{C}^n}
\renewcommand{\C}{\mathbf{C}}
\renewcommand{\Rn}{\mathbf{\R}^n}
\renewcommand{\R}{\mathbf{R}}
\renewcommand{\Pr}{\mathbb{P}}
\renewcommand{\qedsymbol}{$\blacksquare$}
\renewcommand{\tr}{{\rm tr}}
\renewcommand{\vec}{\text{vec}}
\renewcommand{\S}{\mathbf{S}}

\newcommand{\Hn}{\mathbf{H}^n}
\newcommand{\real}{{\rm Re}}
\newcommand{\imag}{{\rm Im}}
\newcommand{\Ck}{\mathbf{C}^k}
\newcommand{\jcomp}{\mathbf{j}}
\newcommand{\Lb}{\mathcal{L}}
\newcommand{\Rk}{\mathbf{R}^k}
\newcommand{\Jcal}{\mathcal{J}}
\newcommand{\Nb}{\mathbf{N}}
\newcommand{\Hb}{\mathbf{H}}
\newcommand{\qeda}{\tag*{$\blacksquare$}}
\newcommand{\SO}{\mathbf{SO}}
\newcommand{\Aps}{A^+}
\newcommand{\Ans}{A^-}
\newcommand{\pmat}[1]{ \begin{pmatrix} #1 \end{pmatrix}}
\newcommand{\Feas}[1]{\mathcal{F}(#1)}
\renewcommand{\Rm}{\mathbb{\R}^m}
\newcommand{\gsf}{\mathsf{g}}
\newcommand{\ol}[1]{\overline{#1}}
\newcommand{\ul}[1]{\underline{#1}}
\newcounter{mytempeqncnt}
\newcommand{\mul}{\mu_{\text{low}}}
\newcommand{\mum}{\mu_{\text{med}}}
\newcommand{\muh}{\mu_{\text{high}}}
\newcommand{\herm}{H}
\newcommand{\Vg}{\mathcal{V}_G}
\newcommand{\VL}{\mathcal{V}_L}
\newcommand{\nG}{n_G}
\newcommand{\nL}{n_L}
\newcommand{\dnom}{d}

\usepackage{mathtools}
\usepackage{empheq} 
\usepackage{cases}
\usepackage{blkarray}

\usepackage{stfloats}
\fnbelowfloat

\usepackage{float}
\usepackage{subfloat}
\usepackage{subcaption}

\usepackage{ifthen}
\newboolean{showcomments}
\setboolean{showcomments}{false}
\newboolean{showcommentsRevision}
\setboolean{showcommentsRevision}{true}

\newcommand{\oldtext}[1]{\ifthenelse{\boolean{showcomments}}
{ \red{ #1}}{}}

\newcommand{\oldtextR}[1]{\ifthenelse{\boolean{showcommentsRevision}}
{ \red{ #1}}{}}

\newcommand{\oldtextF}[1]{\ifthenelse{\boolean{showcommentsRevision}}
{ \red{ #1}}{}}

\begin{document}

\title{\Huge  Robust AC Optimal Power Flow}

\author{Raphael Louca  \qquad Eilyan Bitar
\thanks{Supported in part by NSF grant ECCS-1351621,   NSF grant IIP-1632124,  US DOE under the CERTS initiative, and the Simons Institute for the Theory
of Computing. This paper builds on the authors' preliminary results published as part of the 2016 IEEE Conference on Decision and Control (CDC) \cite{louca2016stochastic}. The current manuscript differs significantly
from the conference version in terms of detailed numerical analyses, formal mathematical proofs, and extended technical discussions.}
\thanks{R. Louca and E. Bitar  are with the School of Electrical and Computer Engineering, Cornell University, Ithaca, NY, 14853, USA.  Emails:  {\tt\small rl553@cornell.edu}, {\tt\small eyb5@cornell.edu}}
}

\maketitle

\begin{abstract} 
There is a growing need for new  optimization methods to facilitate the reliable and cost-effective  operation of power systems  with intermittent renewable energy  resources. In this paper, we formulate the \emph{robust AC optimal power flow} (RAC-OPF) problem as a two-stage robust optimization problem with recourse.
This problem amounts to a nonconvex infinite-dimensional  optimization problem that is computationally intractable, in general. \blueF{Under the assumption that there is adjustable generation or load at every bus in the power transmission network}, we develop a technique to approximate RAC-OPF from within by a finite-dimensional semidefinite program by restricting the space of recourse policies to be \emph{affine} in the uncertain problem data.
We establish a sufficient condition under which the semidefinite program returns an affine recourse policy that is guaranteed to be feasible for the original RAC-OPF problem. 
We illustrate the effectiveness of the proposed optimization method on the \blueF{WSCC 9-bus and IEEE 14-bus test systems}
with different levels of renewable resource penetration and uncertainty. 

\end{abstract}


\section{Introduction} \label{sec_introduction}
\blueN{The AC optimal power flow (AC-OPF) problem is a fundamental
decision problem that is at the heart of power system operations \cite{carpentier1962contribution}.  The AC-OPF problem  is a nonconvex optimization problem,
where the objective is to minimize the cost of generation
subject to power balance constraints described by Kirchhoff’s
current and voltage laws, and operational constraints reflecting
real and reactive limits on power generation, bus voltage magnitudes, and power flows along transmission lines
flows. It is also common to enforce contingency constraints
to ensure that the power system can withstand sudden disturbances,
such as generator or line outages  \cite{stott2012optimal}. 
The nonconvexity of the
AC-OPF  problem is in part due  to the need to enforce
quadratic constraints, which are indefinite in the vector of complex bus
voltages. The treatment of such nonconvexities in the AC-OPF
problem has traditionally relied on the use of local constrained-optimization methods, or the use of approximate linear
models of power flow (e.g., DC-OPF) to convexify the feasible region of the underlying optimization
problem\cite{taylor2015convex}. More recently, considerable
effort has been made to identify conditions under which an
optimal solution to AC-OPF  can be obtained from a 
solution to its semidefinite programming relaxation \cite{lavaei2012zero, low2014convexI, low2014convexII, taylor2015convex}.
According to the US 
Federal Energy Regulatory Commission, a  5\% increase in the efficiency of algorithms
for AC-OPF will yield six billion dollars in savings per year in the United States alone \cite{castillo2013optimal}.
}

Recently, growing concerns about climate change have led many US states to implement policies, which require that a large fraction of their electricity come from  renewable energy resources such as wind and solar. A fundamental  challenge facing the deep integration of such resources stems from the need to  accommodate the intrinsic uncertainty in their power supply. Doing so efficiently will require the development of robust optimization methods for the AC optimal power flow (AC-OPF) problem.
In its most elemental formulation, the \emph{robust AC optimal power flow} (RAC-OPF) problem gives rise to a two-stage robust optimization problem, in which the system operator must determine a day-ahead generation schedule that minimizes the expected cost of dispatch, given an opportunity for recourse to adjust its day-ahead schedule in real-time when the uncertain system variables have been realized (e.g., the power that can be supplied from wind and solar resources).
The RAC-OPF problem  is a nonconvex, infinite-dimensional optimization problem in its most general form due to the nonconvexity of the underlying AC power flow constraints, and the need to optimize over an infinite-dimensional recourse policy space.

\emph{Related Work:} \ In order to treat the nonconvexity of the RAC-OPF problem,  a large fraction of the literature prescribes techniques that rely on a DC linear approximation of the power flow model \cite{bienstock2014chance, jabr2013adjustable, jabr2015robust, lin2018decentralized, lorca2015adaptive, munoz2014piecewise, roald2013analytical, roald2015optimal, summers2014stochastic, warrington2013policy}. \blueN{The computational intractability associated with the need to optimize over an  infinite-dimensional recourse policy space is predominantly addressed by employing affine or piecewise-affine  approximations  of the recourse policy space \cite{bai2015robust, bienstock2014chance,  jabr2013adjustable, jabr2015robust, lin2018decentralized, lorca2015adaptive,  munoz2014piecewise, roald2013analytical, roald2015optimal, summers2014stochastic, vrakopoulou2013probabilistic,  warrington2013policy}.}
Lastly, both \textit{robust}  \cite{bai2015robust, jabr2013adjustable, jabr2015robust, munoz2014piecewise, phan2014two, warrington2013policy} and \textit{chance-constrained} \cite{bienstock2014chance, summers2014stochastic, roald2017chance, roald2013analytical, roald2015optimal, vrakopoulou2013probabilistic} formulations have been proposed to treat the uncertainty in the constraints, which define the RAC-OPF problem.
As the title of this paper suggests, we adopt a robust approach to constraint satisfaction.

The chance-constrained paradigm assigns a distribution to
the uncertain variables, and the system constraints are enforced up
to a prespecified  probability level. The procurement of such distributions, however, is challenging
as  distributions describing  uncertain power system parameters, such
as renewable energy generation, can be difficult to identify \cite{hamon2013value}.
In order to account for the potential inaccuracy in
the specification of the underlying distribution, the authors of \cite{zhang2017two} consider \textit{ambiguous chance constraint formulations} using DC power flow models. Under this paradigm, the underlying distribution is assumed to belong to a closed ball centered around a known distribution. An alternative  treatment of ambiguous chance constraints is provided in \cite{tseng2016random}.
In addition to the basic difficulty inherent to the identification of an accurate distributional model,  chance constraints do
not explicitly account for the magnitude of constraint violations when they occur---however small their a priori specified probability.
To account for this, the conditional value at risk (CVaR) is commonly used in the place of chance constraints
to promote solutions that  minimize the expected magnitude of such constraint violations  \cite{summers2014stochastic, roald2015optimal}.
In contrast to
the chance-constrained paradigm, robust optimization takes a
deterministic approach to uncertainty modeling. In particular,
the uncertain parameters are assumed to vary in a known
and bounded uncertainty set. A robust solution is one that
optimizes the objective function and remains feasible for any
realization of the uncertain parameters in the given set. This immunity
of robust solutions, however, comes at the expense of potential conservatism in their performance.
Under certain mild assumptions, robust linear constraints arising in approximate DC power flow models admit equivalent reformulations as finite-dimensional conic linear constraints  \cite{jabr2013adjustable}. In general, however, robust nonlinear constraints, which arise in AC power flow models, do not admit such equivalent reformulations. The predominant approach to their treatment relies on scenario or sample-based approximation techniques \cite{capitanescu2012cautious}, which give rise to outer approximations (relaxations) of the robustly feasible set.

\emph{Summary of Results:} \  In this paper, we formulate RAC-OPF as a two-stage robust optimization problem with recourse. The formulation considered  departs  from the majority of the extant literature given its treatment of the full AC power flow model. To our knowledge,  the only other papers in the literature that treat the AC model in the robust optimization framework are \cite{bai2015robust, lorca2018adaptive, vrakopoulou2013probabilistic, phan2014two}. They address the nonconvexity, which arises from the AC power flow equations, by means of a  convex (second-order cone or semidefinite) relaxation.\footnote{We refer the reader to \cite{ low2014convexI, low2014convexII,taylor2015convex} for a detailed derivation and analysis of semidefinite and second-order cone relaxations for the \emph{deterministic} (zero-recourse) AC-OPF problem.}  A crucial assumption made in these papers is  the exactness of their convex relaxations.  Exactness of such relaxations for the  RAC-OPF problem  is not guaranteed, and, in particular,  the solutions generated by these relaxations are not guaranteed to be feasible for RAC-OPF.
In this paper,  we adopt an approach that relies on the  restriction of the space of recourse policies to those which are \emph{affine} in the uncertain problem data. Under this restriction and \blueF{the assumption that there is adjustable generation or load at every bus in the power transmission network}, we develop a technique to approximate the  RAC-OPF  problem from within by a finite-dimensional semidefinite program. We establish a sufficient condition under which the resulting semidefinite program---a convex inner approximation to RAC-OPF---yields recourse policies, which are  \emph{guaranteed to be  feasible} for   RAC-OPF. 

\emph{Organization:} \ The paper is organized as follows. In Sections \ref{sec_preliminaries} and \ref{sec:RACOPF}, we  develop the power system model and provide a detailed formulation of RAC-OPF, respectively. In Section  \ref{sec_main_results}, we offer a detailed derivation of the semidefinite programming inner approximation of RAC-OPF, and provide a sufficient condition under which the resulting approximation is guaranteed to have a nonempty feasible region. In Section \ref{sec:mainresults_III}, we describe an iterative optimization method that generates  a  sequence of feasible affine recourse policies with nonincreasing costs. Finally, we  illustrate the effectiveness of the proposed optimization method on \blueF{the WSCC 9-bus and IEEE 14-bus power systems} with different levels of renewable resource penetration and uncertainty. Section \ref{sec:conc} concludes the paper. 
Detailed proofs for the theoretical results
contained in this paper can be found in the Appendix.

\emph{Notation:} \ 
Let $\R$, $\C$, and $\Nb$ denote the set of real, complex, and natural numbers, respectively. 
Given $m \in \Nb$, define the set of natural numbers from 1 to $m$ by $[m]\coloneqq \{1,2,\dots,m\}$. Let $e_i$ denote the real $i^{\rm th}$ standard basis vector, whose dimension will be apparent from the context. Given any pair of complex numbers $z_1,z_2 \in \C$, we write $z_1 \le z_2$ if and only if $\real\{z_1\} \le \real\{z_2\} \text{ and } \imag\{z_1\} \le \imag\{z_2\}$. Given a matrix $X \in \C^{m\times n}$,  we denote its conjugate transpose by $X^*$, and  its $(i,j)$ entry by
$[X]_{ij}$. \blueN{Given a matrix $X \in \R^{m \times n}$, we denote its transpose by $X^T$ and  use $X \ge 0$ to mean that $X$ is entrywise nonnegative. } We denote the trace of a matrix $X$ by $\tr(X)$.  Denote by $I_n$ the $n \times n$ identity matrix.
Denote by $\Hn$ the set of $n\times n$ Hermitian matrices.  We use  $X \succeq 0$ to mean that the matrix $X \in \Hn$ is positive semidefinite.
Finally, let  $\Lb_{k,n}^2$ denote the space of all square-integrable, Borel measurable functions from $\Rk$ to $\Cn.$ A complex-valued function $f$ on $\Rk$ is said to be Borel measurable if both $\real\{f\}$ and $\imag\{f\}$ are real-valued Borel measurable. \blueN{We summarize frequently used symbols and variables in Appendix \ref{app:table_notation}.}

\section{Power System Model} \label{sec_preliminaries}
\blue{In this Section, we develop the \emph{robust AC optimal power flow} (RAC-OPF) problem. We consider a power system consisting of a heterogeneous mix of generators, which differ in terms of their predictability and controllability. We adopt the perspective of the independent system operator (ISO), whose aim is to dispatch generators in order to minimize the expected cost of  serving demand, while respecting generation and transmission capacity constraints---an optimization problem that belongs to the  class of \emph{security-constrained optimal power flow} problems.
 We consider an optimization model that consists of two stages: day-ahead (DA) and real-time (RT). In the day-ahead stage, the ISO determines an initial dispatch of generators subject to uncertainty in certain system variables, e.g., the  supply that will be available from renewable energy resources in real-time.
\blueN{Day-ahead scheduling decisions are critical, as certain generators (e.g., nuclear and coal) have  ramping constraints that limit the extent to which they can adjust their power output in real-time. Accordingly, such ramp-constrained generators must be scheduled to produce well in advance of the delivery time, and, therefore, prior to the realization of certain a priori uncertain variables, e.g., wind power availability.  In the real-time stage, all a priori uncertain variables are realized, and the ISO is given a recourse opportunity to adjust its DA schedule to balance the system at minimum cost. The ramping constraints dictate the extent to which each generator can  adjust its power output around its day-ahead set-point in real-time.} Essentially, the calculation of  a DA schedule that  minimizes the expected cost of generation---subject to  optimal recourse  in real-time---entails the solution of a robust optimization problem with recourse. We provide a precise formulation of this problem in  Section \ref{sec:RACOPF}. }

\subsection{AC Power Flow Model}
We consider a  power transmission network described by an undirected graph $\Gcal \coloneqq (\Vcal,\Ecal)$. The set  of vertices $\Vcal := [n]$  index the transmission  buses, and the set 
of edges $\Ecal \subseteq \Vcal \times \Vcal$  index the transmission lines between buses. We assume that $(j,i) \in \Ecal$ if and only if $(i,j) \in \Ecal$.

The \emph{AC power flow} equations are formulated according to Kirchhoff's voltage and
current laws, which relate bus power injections to voltages \cite{arthur2000power}. We denote by  $Y \in \C^{n \times n}$ the bus admittance matrix,  $s \in \Cn$  the vector of complex (net) bus power injections (generation minus demand), and $v \in \Cn$  the vector of complex bus voltages. The AC power balance equations are given by
\begin{align} \label{eq:power_balance}
		 s_i    =  v^\herm S_iv,
\end{align}
 where $S_i \coloneqq Y^\herm e_ie_i^T$ for all $ i \in \Vcal$.  For each line $(i,j) \in \Ecal$, we denote by $s_{ij} \in \C$ the complex power flow from bus $i$ to  bus $j$. It satisfies
\begin{equation}\label{eq:power_flow}
\begin{alignedat}{4}
	s_{ij}  =   v^\herm S_{ij}v,
\end{alignedat}
\end{equation}
where $S_{ij} \coloneqq e_i e_i^T(\widehat{y}_{ij}/2-[Y]_{ij})^\herm + e_je_i^T[Y]_{ij}^\herm$ 
for all $(i,j) \in \Ecal$.  Here,  $\widehat{y}_{ij} \in \C$ represents the total shunt admittance of line $(i,j)$. 

We require that the following constraints be enforced. The first set of constraints  limit  the range of acceptable voltage magnitudes at each bus $i \in \Vcal$, 
\begin{align*}
	v_i^{\min}  \le |v_i| \le \  v_i^{\max},
\end{align*}
where  $v_i^{\min}\in \R$ and $v_i^{\max} \in \R$ denote  lower and upper limits, respectively, on the bus $i$ voltage magnitude. We also consider active power flow capacity constraints on the transmission lines.  Namely, for each transmission line $(i,j) \in \Ecal$,  the active power flow from bus $i$ to bus $j$ is required to satisfy
\begin{align*}
- \ell_{ij}^{\max} \le v^\herm P_{ij} v  \le   \ell_{ij}^{\max}.
\end{align*}
Here,  we define $P_{ij} \coloneqq (S_{ij}+S_{ij}^\herm)/2$, and let  $\ell_{ij}^{\max} \in \R$ denote the active power flow capacity of transmission line $(i,j)$. \blueN{We note that---within the framework of this paper---it is also possible to accommodate transmission line constraints that limit the magnitude of the current flowing through a line,  as these also give rise to quadratic inequality constraints in the vector of complex bus voltages. The treatment of apparent power flow constraints, however, is beyond  the scope of  the present paper, as these amount to quartic inequality constraints in the vector of complex bus voltages.
}

\subsection{Uncertainty Model}\label{sec:II_Uncertainty}
 All of the  `uncertain' quantities appearing in this paper are described according to the random vector $\bm{\xi}$, which is defined according to the probability space $(\Rk,\Bcal(\Rk),\Pr)$. Here, the  Borel $\sigma$-algebra $\Bcal(\Rk)$  is the set of all events that are assigned probabilities by the  measure $\Pr$. We denote the first and second-order moments of $\bm{\xi}$ by 
 \begin{align*}
 \mu := \Exp[\bm{\xi}]  \ \text{and} \ M := \Exp[\bm{\xi}\bm{\xi}^T],
 \end{align*}
 where  $\Exp[\cdot]$ denotes the expectation operator with respect to $\Pr$.  
Adopting a standard notational convention, we will use  $\xi$ (normal face) to denote  realizations taken by the random vector $\bm{\xi}$ (bold face). We assume throughout the paper  that   the support of the  random vector $\bm{\xi}$  is nonempty, compact, and representable as
\begin{align} \label{eq:supportset}
	\Xi \coloneqq \{ \xi \in \Rk \ | \  \xi_1 =1, \ \xi^TW_j\xi \ge 0, \ j\in [\ell]\}.
\end{align}
 Here, each matrix $W_j \in \R^{k\times k}$ is defined according to
\begin{align} \label{eq:Wj_repre}
	W_j \coloneqq \bmat{\omega_j & w_j^T \\[0.2em] w_j & -\Omega_j^T\Omega_j},
\end{align}
where $\omega_j \in \R$, $w_j \in \R^{k-1}$, and $\Omega_j \in \R^{n_j\times(k-1)}$  for some $n_j \in \Nb$. It is important to note that the representation of the support in \eqref{eq:supportset} is general enough to describe any subset of the hyperplane $\{ \xi \in \Rk \ | \  \xi_1 =1 \}$  that is defined according to a finite intersection of
arbitrary half spaces  and ellipsoids.  We will occasionally refer to the support set  $\Xi$  as the `uncertainty set' associated with the random vector $\bm{\xi}$. \blueN{Naturally, the particular specification of the uncertainty set  will  depend on the specific nature of the  uncertain quantities that   the random vector $\bm{\xi}$ is meant to capture, e.g., the maximum active power that  a wind generator is able to produce at a node.}

Some remarks regarding our uncertainty model are in order. First, the requirement  that $\xi_1=1$ for all $\xi \in \Xi$ is for notational convenience, as it allows one to represent affine functions of $(\xi_2,\dots,\xi_k)$ as linear functions of $\xi$. Second,  it is important to emphasize that all of the results contained in this paper  depend  on the probability distribution of the random vector $\bm{\xi}$ only through its support, mean, and second-order moment. No additional information about the distribution is required. 
The following mild technical assumption  is assumed to hold throughout the paper.
\begin{assumptio} \label{ass:slater}\blueN{There exists  $\xi \in \Xi$} such that $\xi^T W_j\xi > 0$ for all $j \in [\ell]$.
\end{assumptio}
The assumption that the support set $\Xi$ admits a strictly feasible point will prove useful to the derivation of our subsequent theoretical results, as it ensures that  $\Xi$ spans all of $\R^k$. This, in turn, guarantees that the second-order moment matrix $M$ is positive definite and invertible. We refer the reader to \cite[Prop. 2]{kuhn2011primal} for a proof of this claim.

\subsection{Generator and Load Models}
\blue{\emph{Load Model:}  \  The real-time  demand for power at each bus $i\in \Vcal$ is assumed to be
fixed and known. We  denote it by $\dnom_i  \in \C$. 
}

\vspace*{.08in}

\blue{\emph{Generator Model:} \ For ease of exposition,  we assume throughout the paper that there is at most a single generator at each bus $i \in \Vcal$. We consider a generator model in which the real-time (RT) supply of power, as determined by the ISO, is allowed to depend on the realization of the random vector $\bm{\xi}$. 
Hence, we denote the power produced  by generator $i$ in real-time by  $g_i(\bm{\xi})$, where  $g_i \in \Lb^2_{k,1}$ is a \emph{recourse function} determined by the ISO for each generator $i$. The cost incurred by each generator $i$ for producing $g_i(\bm{\xi})$ is assumed to be  linear in the active power produced. It is defined as
\begin{align*}
 \alpha_i \real\{g_i(\bm{\xi})\}, \qquad i \in \Vcal,
\end{align*}
where $\alpha_i \ge 0$  denotes the  marginal cost of  real power generation at bus  $i \in \Vcal$.}

\blue{We consider a generator model in which the power capacity available to each generator in real-time is allowed to be uncertain day-ahead. Namely, the power produced by each generator $i \in\Vcal$ in real-time must satisfy
\begin{align} \label{eq:RTcap}
 \underline{g}_i(\bm{\xi})  \leq  g_i(\bm{\xi}) \leq \overline{g}_i(\bm{\xi}) , \qquad i \in \Vcal.
\end{align}
Here,  $\underline{g}_i \in \Lb^2_{k,1}$ and  $\overline{g}_i \in \Lb^2_{k,1}$ denote  lower and upper bounds on the power produced by  generator $i$ in real-time. Such a model is general enough to describe apriori uncertainty in renewable power supply, as well as unscheduled generator outages. We assume that the random generation capacities satisfy
$$g_i^{\min} \le  \underline{g}_{i}(\bm{\xi})  \le  \overline{g}_i(\bm{\xi})  \le  g_i^{\max}, \qquad  i \in \Vcal.$$
Here, $g_i^{\min} \in \C$ and  $\ g_i^{\max} \in \C$ denote the \emph{nameplate} minimum and maximum capacities of generator $i$, respectively. The corresponding vectors are denoted by  $\underline{g}(\bm{\xi})$,  $\overline{g}(\bm{\xi})$,   $g^{\min}$,   and   $g^{\max}$. 
}

In practice,  generators have ramping constraints that limit the extent to which they can adjust their power production in real-time. 
We model the limited ramping capability of each generator $i$ according to the following  pair of constraints
\begin{align} \label{eq:ramp}
		r^{\min}_i  \le  g_i(\bm{\xi})-g_i^0  \le  r^{\max}_i , \qquad i \in \Vcal,
\end{align}
where  $r^{\min}_i \in \C $ and  $r^{\max}_i \in \C$ represent the ramp-down and ramp-up limits, respectively, associated with generator $i$. Here, 
$g_i^0 \in \C$ denotes generator $i$'s   day-ahead (DA) dispatch, also determined by the ISO. The DA dispatch of each generator is required to satisfy its nameplate generation capacity constraints given by
\begin{align} \label{eq:DAcap}
		g_i^{\min} \le g_i^0  \le \ g^{\max}_i , \qquad i \in \Vcal.
\end{align}

\begin{exampl}[Generator types]   \em
The generator model that we consider in this paper captures a wide range of generator types. We provide several important examples below. Let $g_i^0$ be a DA dispatch level satisfying \eqref{eq:DAcap}.  Generator $i$ is said to be:
\begin{itemize}
\item \textit{Completely inflexible} (e.g., nuclear) if  its  RT power output is restricted to $$g_i(\bm{\xi})  = g_i^0.$$
\item \textit{Completely flexible} (e.g.,  gas, oil) if  its RT  power output is restricted to $$g_i^{\min} \le  g_i(\bm{\xi})  \le  g_i^{\max}.$$
\item \textit{Intermittent} (e.g., wind, solar) if  its RT power output is restricted to $$\underline{g}_i(\bm{\xi}) \le  g_i(\bm{\xi})  \le  \overline{g}_i(\bm{\xi}).$$
\end{itemize}
\end{exampl}

We state the following technical assumption, which requires that the  RT generation capacities exhibit a linear dependence on the random vector $\bm{\xi}$. Assumption \ref{ass:lin} is required to hold for the remainder of the paper.

\begin{assumptio} \label{ass:lin} There exist matrices  $\underline{G}  \in \C^{n \times k}$ and $\overline{G} \in \C^{n \times k}$ such that  $\underline{g}(\bm{\xi}) = \underline{G} \bm{\xi}$ and  $\overline{g}(\bm{\xi}) = \overline{G}\bm{\xi}$.
\end{assumptio}

\section{Robust AC Optimal Power Flow} \label{sec:RACOPF}
Building on the previously defined models, we formulate the \emph{robust AC optimal power flow} (RAC-OPF) problem as follows.
	\begin{alignat}{4}\label{prog:SDR_OPF_polar}
		& \underset{}{\text{minimize}}   
		&& \Exp \left[  \sum_{i=1}^n \alpha_i {\rm Re}\{g_i(\bm{\xi})\}  \right]   \\
		& \text{subject to}  \ \    
		&&  g^0 \in \C^n, \ \  g \in \Lb^2_{k,n}, \ \  v \in \Lb^2_{k,n}     \nonumber \\[0.5em]
		&&& g_i^{\min}  \le  g_i^0 \le  g_i^{\max}, &&    i\in  \Vcal   \nonumber \\[0.5em]
		&&& \underline{g}_i(\xi)   \le  g_i(\xi)   \le  \overline{g}_i(\xi),   &&    i\in  \Vcal 
				\makebox[42pt][r] {\smash{\raisebox{-2.7\baselineskip}{$\left.\rule{0pt}{3.335\baselineskip}\right\}\hspace{-0.05cm}\forall \ \hspace{-0.05cm}\xi \hspace{-0.05cm}\in \hspace{-0.05cm} \Xi. $}}}
		 \nonumber \\[0.5em]
		&&&  r_i^{\min} \le  g_i(\xi) - g^0_i  \le  r_i^{\max}, &&   i\in  \Vcal   \nonumber \\[0.5em]
		&&&  g_i(\xi) -v(\xi)^HS_iv(\xi)  =  \dnom_i, &&   i\in  \Vcal   \hspace{.2in}  \nonumber \\[0.5em]
		&&&  v_i^{\min}  \le  |v_i(\xi)|  \le  v_i^{\max},  &&  i\in  \Vcal  \nonumber \\[0.5em]
		&&&  |v(\xi)^H P_{ij} v(\xi)|  \le   \ell_{ij}^{\max}, &&   \hspace{-.22in} (i,j)\in \Ecal   \nonumber
\end{alignat} 
As  previously described, the RAC-OPF problem amounts to a two-stage robust  optimization problem with recourse. The single-period formulation of RAC-OPF  that we consider is similar in structure to the single-period formulations studied in \cite{bai2015robust, jabr2013adjustable, pritchard2010single, vrakopoulou2013probabilistic}.
We briefly summarize the timing and structure of the  decision variables  and constraints of  the RAC-OPF problem.
\begin{itemize} \itemsep0.6em
\item The \emph{first-stage} (day-ahead) decisions entail  the determination of  a DA generator dispatch  $g_0 \in \C^n$ subject to optimal recourse in the second stage, which will adjust the DA dispatch given a   realization of the random vector $\bm{\xi}$.

\blue{
\item In the \emph{second-stage} (real-time), the random vector $\bm{\xi}$ is realized, and the ISO is given
a recourse opportunity to adjust its DA generator dispatch to balance the system at minimum cost. The second-stage decision entails the determination of the  RT generator dispatch \blueN{$g \in \Lb^2_{k,n}$}  and the RT bus voltages  \blueN{$v \in \Lb^2_{k,n}$}. 
}

\item All decisions must be jointly determined in such a manner as to (i)  minimize the expected cost of generation, and (ii) guarantee that all system constraints are satisfied given \emph{any} realization $\xi \in \Xi$ of the random vector $\bm{\xi}$ in real-time, i.e., robust constraint satisfaction.
\end{itemize}

\blueN{
\begin{remar}[Minimax formulation] \label{rem:mini_max} While the formulation of the RAC-OPF problem in \eqref{prog:SDR_OPF_polar} entails minimizing the expected cost of generation, the computational methods and theoretical results developed in Sections \ref{sec_main_results} -\ref{sec:mainresults_III} can be  generalized to accommodate a minimax formulation of the RAC-OPF problem, which entails minimizing the maximum (worst-case) cost of generation:
\begin{align} \label{eq:minimax}
	\max_{\xi \in \Xi} \ \Bigg\{ \sum_{i=1}^n \alpha_i {\rm Re}\{g_i(\xi)\} \Bigg\}.
\end{align}
The resulting minimax formulation of the RAC-OPF problem under the objective function \eqref{eq:minimax} can be equivalently reformulated to resemble  \eqref{prog:SDR_OPF_polar} by putting it in its \emph{epigraph form}. 
\end{remar}
}

\subsection{Concise Formulation of RAC-OPF} 
It will be convenient to our analysis in the sequel to work with a more concise representation of the RAC-OPF problem.  We do so by first eliminating the RT generator dispatch  variables $g\in \Lb^2_{k,n}$  through their direct substitution according to the nodal power balance  equations. Second, by redefining the DA generator dispatch $g_0 \in \C^n $  as a real vector $x := [ \real\{g_0\}^T, \imag\{g_0\}^T]^T $, one can rewrite problem \eqref{prog:SDR_OPF_polar} more compactly in the following form:
\begin{alignat}{4} \label{prog:two_stage_general}
	& \underset{}{\text{minimize}} \ \ \
	& & \Exp \left[ v(\bm{\xi})^H A_0 v(\bm{\xi})\right] \tag{$\Pcal$}  \\
	& \text{subject to}     
	&&  x \in \R^{2n}, \ v \in \Lb^2_{k,n} \nonumber  \\[0.3em]
	&&&  v(\xi)^H A_iv(\xi) +  b_i^T x \le c_i^T \xi,  \ \  i \in [m], \ \ \forall \  \xi  \in  \Xi \nonumber \\[0.3em]
	&&& Ex  \le f, \nonumber 
\end{alignat} 
where $m:=10n+2|\Ecal|$.  \blueN{We remark that in the above reformulation of problem \eqref{prog:SDR_OPF_polar}, we have eliminated the constant term $\sum_{i=1}^n\alpha_i\real\{\dnom_i\}$ from the objective function as this does not affect the optimal solution of the RAC-OPF problem.}
 It is straightforward to construct the matrices  $E \in  \R^{4n\times 2n}$,  $f \in \R^{4n}$,   $A_i \in \Hn$ ($i = 0,\dots, m$),     $b_i \in \R^{2n}$  ($i = 1,\dots, m$),  and $c_i \in \R^k$  ($i = 1,\dots, m$) given the underlying  problem data  specified in the  RAC-OPF  problem \eqref{prog:SDR_OPF_polar}.  We refer the reader to Appendix  \ref{app:matricesAi} for their specification. 

\begin{remar}[Eliminating quadratic equality constraints] 
We remark that the formulation of the original  RAC-OPF problem \eqref{prog:SDR_OPF_polar}  assumes that there is adjustable generation at every bus in the power transmission network. The advantage of this \blueF{rather limiting assumption} is that it enables the elimination of all  nodal power
balance equality constraints  in the equivalent reformulation of the RAC-OPF problem given by \ref{prog:two_stage_general}. The ability to eliminate these nonconvex quadratic equality constraints will be  essential to the convex inner approximation technique developed in Section \ref{sec:main_results_convexification}. In Appendix   \ref{app:Load_Shedding}, we provide an alternative formulation of the RAC-OPF problem to accommodate the treatment of more general power systems in which \emph{load shedding} is permitted at non-generator (load) buses, where any reduction in load is penalized according to the value of lost load (VOLL)). \blueF{We refer the reader to Section \ref{sec_numerical_studies} for several numerical case studies, which asses the extent to which the allowance of load shedding manifests in an actual reduction in load under the dispatch policies proposed in this paper. 
Finally, we note that the approximation technique proposed in this paper cannot be applied to power systems with transmission buses that have neither adjustable generation or load.}
\end{remar}

\section{Convex Inner Approximation of RAC-OPF}
\label{sec_main_results}

Problem \ref{prog:two_stage_general}  is computationally intractable, in general, as it is both \emph{infinite-dimensional}  and \emph{nonconvex}. The nonconvexity is due, in part, to the feasible set,  which is defined by a number of indefinite quadratic  inequality constraints in the vector of complex bus voltages. 
The infinite-dimensionality of the optimization problem \ref{prog:two_stage_general}  derives from both the infinite-dimensionality of the recourse decision variables, and the infinite number of constraints due to the infinite cardinality of the uncertainty set $\Xi$. In what follows,  we develop a systematic approach to approximate problem \ref{prog:two_stage_general} from within by a finite-dimensional semidefinite program, and provide a sufficient condition under which the resulting inner approximation is guaranteed to have a nonempty feasible region. The proposed method for approximation centers on the restriction of the infinite-dimensional space of recourse policies to those which are \emph{linear} in the random vector $\bm{\xi}$.  

\subsection{Affine Recourse Policies} \label{sec:main_results_I}
As the initial step in the derivation of  a tractable inner approximation to problem  \ref{prog:two_stage_general}, we first restrict the functional form of the recourse decision variables (i.e., the complex bus voltages) to be  linear in the random vector $\bm{\xi}$.\footnote{We note that this  restriction on the functional form of the complex bus voltages implies a \emph{quadratic} dependency of the nodal power generation levels on the random vector $\bm{\xi}$. That is to say,  $g_i(\bm{\xi})  = \bm{\xi}^T V^H S_i V \bm{\xi}  +  d_i$ for each node $i \in \Vcal$.} That is to say, we require that
\begin{align}\label{eq:affine_decision}
		v(\bm{\xi})   = V \bm{\xi},
	\end{align}
where $V \in \C^{n \times k}$. This restriction to affine recourse policies gives rise to the following optimization problem \ref{prog:two_stage_affine_policies}, which stands as an inner approximation to the original problem  \ref{prog:two_stage_general}.  
\begin{alignat}{4} \label{prog:two_stage_affine_policies}
	& \underset{}{\text{minimize}} \tag{$\Pcal_{\rm I}$}\ \  \
	& & \tr(MV^\herm A_0V) \nonumber \\[0.3em]
	& \text{subject to}     
	&&  x \in \R^{2n}, \ V \in \C^{n \times k} \nonumber \\[0.3em]
	&&&  \xi^T V^\herm A_i V \xi +  b_i^Tx \le c_i^T \xi, \ \  i \in [m], \ \ \forall \  \xi  \in  \Xi \nonumber \\[0.3em]
	&&&   Ex  \le f, \nonumber 
\end{alignat}
We have used linearity of the expectation and trace operators, and the invariance of trace under cyclic permutations to massage the original objective function to obtain $$\Exp[\bm{\xi}^T V^\herm A_0V \bm{\xi}]  =  \Exp[\tr(\bm{\xi}\bm{\xi}^T V^\herm A_0V) ] =   \tr(\Exp[\bm{\xi}\bm{\xi}^T] V^\herm A_0V).$$ 

The resulting problem \ref{prog:two_stage_affine_policies} amounts to a semi-infinite, nonconvex quadratically constrained quadratic program.\footnote{A semi-infinite program is an optimization problem with an infinite number of constraints, and  finitely many decision variables.} More specifically, the \blue{restriction to affine recourse policies yields an optimization problem  that has finite-dimensional decision
variables.  However, Problem \ref{prog:two_stage_affine_policies} remains to be computationally intractable, as it requires the satisfaction of infinitely many constraints due to the continuous structure of the uncertainty set $\Xi$. We address this issue in Lemma \ref{lemma:Sprocedure} by employing weak
duality to obtain a \emph{sufficient} set of finitely many constraints. Such  an approximation of the infinite constraint set can also be derived through a direct application of the so-called $S$-procedure \cite{boyd1994linear}.
 We state Lemma \ref{lemma:Sprocedure} without proof, as it follows directly from Proposition 6 in \cite{kuhn2011primal}.}

\oldtext{ restriction to affine recourse policies results in an optimization problem \ref{prog:two_stage_affine_policies} whose decision variables range over finite-dimensional spaces. However, due to the continuous structure of the uncertainty set $\Xi$,  problem \ref{prog:two_stage_affine_policies} has infinitely many constraints and is, in general, intractable.
To account for this, we employ weak
duality to obtain a \emph{sufficient} set of finitely many constraints. We remark that such  an approximation of the infinite constraint set can also be derived through a direct application of the so-called $S$-procedure \cite{boyd1994linear}.
We have the following result, which follows from Proposition 6 in \cite{kuhn2011primal}.}

\begin{lemm}\label{lemma:Sprocedure}
Let  $P \in \Hb^{k}, \ q \in \R^k$,  $r \in \R$, and  $Q \coloneqq (e_1q^T + qe_1^T)/2$. Consider the following two statements:\vspace{0.1cm}
\begin{enumerate}[(i)] \itemsep0.4em
	\item\label{lemma:S2} $\xi^T P\xi + q^T \xi +r \le 0$ for all $\xi \in \Xi$,

	\item \label{lemma:S1} $\exists$ $\lambda \in \R^{\ell}$ with $\lambda \le 0$ and 
	$
		P+Q+re_1e_1^T -\sum\limits_{j=1}^\ell \lambda_j W_j   \preceq  0,
	$
\end{enumerate}
where $W_j$ is as defined in \eqref{eq:Wj_repre}.
For any $\ell \in \Nb$, it holds that \eqref{lemma:S1} implies \eqref{lemma:S2}. If $\ell=1$, then  \eqref{lemma:S2} and \eqref{lemma:S1} are equivalent.
\end{lemm}
Using Lemma \ref{lemma:Sprocedure}, one can approximate the infinite constraint set of problem \ref{prog:two_stage_affine_policies} from within by finitely many matrix inequality constraints. More precisely, a direct application of Lemma \ref{lemma:Sprocedure} to each of the quadratic constraints in problem \ref{prog:two_stage_affine_policies} gives rise to the following finite-dimensional optimization problem:
\begin{alignat}{8} \label{prog:two_stage_Sprocedure1}
	& \underset{}{\text{minimize}}   
	& & \tr(MV^\herm A_0V) \tag{$\Pcal_{\rm II}$} \\
	& \text{subject to} \ \ \ 
	&&  x \in \R^{2n}, \ V \in \C^{n \times k}, \ \Lambda \in \R^{m\times\ell} \nonumber  \\[-0.3em]
	&&&  V^\herm A_iV - C_i + (b_i^T x)e_1e_1^T - \sum_{j=1}^\ell [\Lambda]_{ij} W_j &&\preceq 0,  \nonumber\\[-0.5em]
	&&& \forall \   i \in [m], \nonumber \\
	&&&   Ex  \le f, \nonumber \\
		&&& \Lambda   \le 0, \nonumber 
\end{alignat} 
where we define $C_i \coloneqq  (e_1c_i^T + c_ie_1^T)/2$ for each $i \in [m]$.

\begin{remar}  \blueF{We remark that $\Lambda \in \mathbf{R}^{m \times \ell}$ is a decision variable in problem $\Pcal_{\rm II}$, which emerges from the from the application of Lemma \ref{lemma:Sprocedure}  to the robust inequality constraints in problem $\Pcal_{I}.$} It follows  from Lemma \ref{lemma:Sprocedure} that problem \ref{prog:two_stage_Sprocedure1} is an inner approximation to problem \ref{prog:two_stage_affine_policies} in general; and is equivalent to problem \ref{prog:two_stage_affine_policies} when $\ell=1$. 
\end{remar}

\subsection{Convexifying the Inner Approximation}\label{sec:main_results_convexification}
Problem \ref{prog:two_stage_Sprocedure1} is a finite-dimensional inner approximation to the original problem \ref{prog:two_stage_general}. It, however, remains to be  nonconvex, because of the indefinite quadratic functions appearing in both the objective and the inequality constraints. In what follows, we develop a  method to convexify problem \ref{prog:two_stage_Sprocedure1} from within by  replacing each indefinite quadratic function with a \emph{majorizing} convex quadratic function. We state the resulting convex program, which approximates  \ref{prog:two_stage_Sprocedure1} from within, in Proposition \ref{prop:Linearization}.

\oldtext{The proposed method is based on the simple observation that each indefinite quadratic function can be decomposed as a sum of a convex quadratic function and a concave quadratic function. 
We then approximate the concave function from above with its linearization at a point. The sum of this linearization with the convex component of the original function yields a convex global overestimator of the original indefinite quadratic function.}

The proposed method is based on the decomposition of an indefinite quadratic function as the difference of convex functions, i.e., the sum of a convex quadratic function and a concave quadratic function. We construct a convex global overestimator of the original indefinite quadratic function by  
linearizing the concave function at a point.\footnote{In Section \ref{sec:main_resultsII}, we provide a method to select the point around which the linearization is constructed  in order to guarantee that the resulting convex feasible set is nonempty.}
More precisely, for each matrix $A_i$, define the decomposition $$A_i  =    \Aps_i  +  \Ans_i,$$
where $\Aps_i \succeq 0$  and   $\Ans_i  \preceq 0$ denote the positive semidefinite and negative semidefinite parts of $A_i$, respectively. Using this matrix decomposition, define the function  $H_i:\C^{n \times k} \times \C^{n \times k} \rightarrow \Hb^k$ according to
\begin{align*}
	H_i(V,Z)  \coloneqq  V^\herm \Aps_iV  +  Z^\herm \Ans_iV + V^\herm \Ans_iZ  -  Z^\herm \Ans_iZ,
\end{align*}
for each $i \in [m]$. The first term of $H_i$ is the convex component of the original quadratic function $V^\herm A_iV$. The remaining terms represent the linearization of the concave component at a point $Z$.
Consequently, for any matrix $Z$, the function $H_i(V,Z)$ is matrix convex in $V$.\footnote{A function $f:\C^{n \times k} \rightarrow 
\Hb^k$ is said to be matrix convex if for all matrices $X,Y$ and $0\le \theta \le 1$, we have $f(\theta X+(1-\theta)Y) \preceq \theta f(X) + (1-\theta)f(Y)$.}   The following result highlights two important properties of $H_i$. Its proof can be found in Appendix B-A. 

\begin{lemm}\label{lemma:properties_Hi}
	Let $Z \in \C^{n \times k}$. For each $i \in [m]$, it holds that 
	\begin{enumerate}[(i)] \itemsep0.6em
		\item \label{lemma:properties_Hi_a} $V^\herm A_iV  \preceq H_i(V,Z), \quad  \forall \ V \in \C^{n \times k}$,
		\item \label{lemma:properties_Hi_b}$\tr(MV^\herm A_iV)  \le  \tr(MH_i(V,Z)), \quad  \forall \ V \in \C^{n \times k}$. 
	\end{enumerate}
\end{lemm}

Property  \eqref{lemma:properties_Hi_a} provides a way of approximating the nonconvex feasible set of problem  \ref{prog:two_stage_Sprocedure1} from withing by a convex set. Property \eqref{lemma:properties_Hi_b}, on the other hand, provides way of majorizing the nonconvex objective of problem \ref{prog:two_stage_Sprocedure1} with a convex function. In Proposition \ref{prop:Linearization}, we employ these approximations to specify a convex program whose optimal solution is guaranteed to be a feasible solution for the original problem \ref{prog:two_stage_general}. Its proof follows directly from Lemma \ref{lemma:properties_Hi}. We, therefore, omit it for the sake of brevity.

\begin{propositio} \label{prop:Linearization} Let $V_0 \in \C^{n \times k}$, and suppose that $(\overline{x},\overline{V},\overline{\Lambda})$ is an optimal solution for the following convex program:
\begin{alignat}{8} \label{prog:two_stage_Linearization}
	& \underset{}{{\rm minimize}}   
	& & \tr(MH_0(V,V_0))  \tag{$\Pcal_{\rm III}(V_0)$} \\
	& {\rm subject \ to} \ \ \ 
	&&  x \in \R^{p}, \ V \in \C^{n \times k}, \ \Lambda \in \R^{m\times\ell} \nonumber \\[-0.3em]
	&&&  H_i(V,V_0) - C_i + (b_i^T x)e_1e_1^T - \sum_{j=1}^\ell [\Lambda]_{ij} W_j &&\preceq 0, \nonumber\\[-0.5em]
	&&& \forall \ i \in [m],  \nonumber \\
	&&& Ex  \le f, \nonumber \\
		&&&\Lambda   \le 0. \nonumber 
\end{alignat}
Define the function $\overline{v} \in \Lcal^{2}_{k,n}$  according to $\overline{v}(\bm{\xi}) = \overline{V}\bm{\xi}$. 
Then $(\overline{x},\overline{v})$  is a feasible solution for the original problem \ref{prog:two_stage_general}.
\end{propositio}
\oldtext{We note that problem \ref{prog:two_stage_Linearization} can be equivalently reformulated as a semidefinite program using the Schur complement condition for positive semidefiniteness. }

We note that problem \ref{prog:two_stage_Linearization} can be equivalently reformulated as a semidefinite program using Schur's complement formula. We refer the reader to Appendix C 
for the details of this reformulation.

\begin{table*}
\setlength{\tabcolsep}{7.4pt}
\renewcommand{\arraystretch}{1.2} 
\begin{tabular*}{7.15in}{lllcccccccccccccc} 
\toprule 
& & & &  Inflexible Gen. & & \multicolumn{5}{c}{Intermittent Gen.} & &  \multicolumn{2}{c}{Flexible Gen.} & & Load Buses    \\  
\multicolumn{1}{c}{Parameters} & & Units & & 1 & & 2 &  3&  4&  5& 6 &    & 7 & 8 & & 9  \\
\cmidrule{1-1} \cmidrule{3-3} \cmidrule{5-5} \cmidrule{7-11} \cmidrule{13-14} \cmidrule{16-16}
$\alpha_i $ & & \$/\text{MW}  &  & 30 & & 0 & 0 & 0 & 0 & 0   &  &    50 & 50 & & -- \\          
$\real\{d_i\}$ & & MW & & 0 & & 0 & 0 & 100 & 0 & 125 &  &  0 & 0 & & 90   \\
$\imag\{d_i\}$ & & MVAR & & 0 & & 0 & 0 & 35 & 0 & 50 &  &  0 & 0  & & 30 \\
$\real\{g^{\max}_i\}$ & & MW & &  200 &  & 30  & 30  & 30 &30  & 30    & & 250  &  270 & & --  \\
$\imag\{g^{\max}_i\}$ & & MVAR & & 300  &  & $\dagger$ & $\dagger$  & $\dagger$   &$\dagger$  & $\dagger$ & & 300 & 300 & & --   \\
$\real\{g^{\min}_i\}$& & MW & & 10   & & 0 & 0  & 0   & 0  & 0   & &  10 &10 & & --   \\
$\imag\{g^{\min}_i\}$& & MVAR & &  -300 & &  $\dagger$ & $\dagger$  & $\dagger$   &$\dagger$   & $\dagger$    & & -300 &-300 & & --   \\
$ \real\{r_i^{\max}\}=-\real\{r_i^{\min}\}$  & & MW & &0 & & $\infty$  & $\infty$ & $\infty$ & $\infty$ &$\infty$  & & $\infty$ &$\infty$ & & --  \\               
$ \imag\{r_i^{\max}\}=-\imag\{r_i^{\min}\}$  & & MVAR & &0 & & $\infty$  & $\infty$ & $\infty$ & $\infty$ &$\infty$ & & $\infty$ &$\infty$ & & -- \\ 
Bus index (per \cite{sauer1998power})   & & & & 2  & & 4 & 6& 7&8&9 & &1&3 & & 5\\    \bottomrule           
\end{tabular*}
\caption{(WSCC 9-bus system). \   Specification of each generator's location, marginal cost, and constraint parameters. The $\dagger$ symbol indicates that the corresponding value in the table is determined by equation \eqref{eq:q_max}. 
} 
\label{table:1a}
\vspace{-0.1cm}
\end{table*}

\subsection{Guaranteeing Nonemptiness of the Inner Approximation}	\label{sec:main_resultsII} 
In order to convexify the RAC-OPF problem according the method  developed in Section \ref{sec:main_results_convexification}, one has to select a matrix $V_0 \in \C^{n \times k}$, \blue{which yields an inner approximation \ref{prog:two_stage_Linearization} with a \textit{nonempty feasible set}. In what follows, we develop a method to compute  such a matrix.}
\oldtext{that results in a \emph{nonempty  feasible set} for the  inner approximation \ref{prog:two_stage_Linearization}. In this section, we provide a method for computing one such matrix.}
The method we propose entails the calculation of a day-ahead dispatch $g^0 \in \C^n$, which is guaranteed to be feasible for the RAC-OPF problem  without requiring adjustment (recourse) in real-time. In order to do so, we must first characterize the guaranteed range of available power supply at each bus.
For each bus $i \in \Vcal$, this amounts to the specification of upper and lower limits $\gamma_i^{\max} \in \C$ and $\gamma_i^{\min} \in \C$, such that 
\begin{align*}
\ul{g}_i(\xi) \leq  \gamma_i^{\min} \leq  \gamma_i^{\max}  \leq  \ol{g}_i(\xi),  \ \ \forall \ \xi \in \Xi.
\end{align*}
We specify these  limits according to
\begin{equation} \label{eq:gammata}
\begin{aligned}
	&\gamma_i^{\min} = \max_{\xi \in \Xi}   \left( \real\{\ul{g}_i(\xi)\}  \right)+ \mathbf{j}\max_{\xi \in \Xi}  \left(  \imag\{\ul{g}_i(\xi)\} \right), \\
	& \gamma_{i}^{\max} = \min_{\xi \in \Xi} \Big( \real\{\ol{g}_i(\xi)\} \Big) +\mathbf{j} \min_{\xi \in \Xi} \Big( \imag\{\ol{g}_i(\xi)\} \Big).
\end{aligned}
\end{equation}
\blueF{
It is important to note that the limits $\gamma_i^{\min}$ and $\gamma_i^{\max}$  can be efficiently calculated, as the  optimization problems in \eqref{eq:gammata} are \emph{convex quadratically constrained quadratic programs}.  This follows from the assumed linearity of the objective function in $\xi$ (cf. Assumption \ref{ass:lin}), and the definition of $\Xi$ as the finite intersection of ellipsoids and half spaces (cf. Eq.  \eqref{eq:supportset}). }

 Using these conservative generation limits, a day-ahead dispatch that is guaranteed to be feasible for the RAC-OPF problem can  calculated by solving the following (deterministic) \blueF{\emph{zero-recourse  AC-OPF}} problem.
\begin{alignat}{4}\label{prog:SDR_OPF_polar1}
		& \underset{}{\text{minimize}}   \ \  \
		&&  \sum_{i=1}^n \alpha_i v^\herm \left(\frac{S_i+S_i^\herm}{2}  \right) v   \\[0.6em]
		& \text{subject to}  \ \    
		&&  v \in \Cn   \nonumber \\[0.5em]
		&&& \gamma_i^{\min}-\dnom_i  \le v^\herm S_i v  \le  \gamma_i^{\max}- \dnom_i,  \quad &&    i \in \Vcal, \nonumber  \\[0.5em]
		&&&  v_i^{\min}  \le |v_i|  \le  v_i^{\max},\quad &&   i \in  \Vcal, \nonumber \\[0.5em]
		&&&  -\ell_{ij}^{\max} \le  v^\herm P_{ij} v  \le   \ell_{ij}^{\max},  &&    (i,j)\in \Ecal. \nonumber
\end{alignat}

The following result shows  that \emph{any feasible solution} to  the \blueF{zero-recourse AC-OPF} problem \eqref{prog:SDR_OPF_polar1} can be used to construct  a matrix $V_0$ that is guaranteed to  induce a convex inner approximation     
\ref{prog:two_stage_Linearization} of the RAC-OPF problem with a  \emph{nonempty feasible region.}
The proof of Proposition \ref{prop:feasibility} can be found in Appendix B-B. 

\begin{propositio} \label{prop:feasibility} Let $v_0 \in \C^n$ be a feasible solution to \eqref{prog:SDR_OPF_polar1}, and define a matrix $V_0 \coloneqq v_0 e_1^T$. It follows that the optimization problem \ref{prog:two_stage_Linearization} has a nonempty feasible region.
\end{propositio}

Several comments are in order.  First, it is important to note that, despite being deterministic, the  \blueF{zero-recourse} AC-OPF problem \eqref{prog:SDR_OPF_polar1} is nonconvex and   computationally intractable, in general. However, there are many off-the-shelf optimization routines (e.g., Matpower \cite{zimmerman2011matpower}) that are effective in producing feasible solutions to problem  \eqref{prog:SDR_OPF_polar1}.  
\blueF{Second, it is also crucial to emphasize that Proposition \ref{prop:feasibility} is only useful if the zero-recourse AC-OPF  problem \eqref{prog:SDR_OPF_polar1} is in fact feasible.  A necessary condition for the feasibility of problem \eqref{prog:SDR_OPF_polar1} is that 
$\gamma_i^{\min} \leq  \gamma_i^{\max}$ for every generator $i \in \Vcal$. This condition  is clearly satisfied by  conventional generators with `firm' real-time generation capacities, i.e., $\ul{g}_i(\xi) = g_i^{\min}$ and   $\ol{g}_i(\xi) = g_i^{\max}$ for all $\xi \in \Xi$. This condition is also satisfied by intermittent generators whose real-time active power supply can be fully `curtailed', i.e., $\real\{\ul{g}_i(\xi)\} = 0$ for all $\xi \in \Xi$. Such an   assumption of curtailable supply is reasonable for modern-day wind and solar power facilities. We refer the reader to Section \ref{sec:rengenmodel} for a specific example of curtailable intermittent generation. It is also important to note that the assumption that intermittent generators be curtailable is necessary, as problem \eqref{prog:SDR_OPF_polar1} may become infeasible for power systems with intermittent generators whose real-time supply  is uncertain in day-ahead and cannot be curtailed in real-time, i.e.,  $\ul{g}_i(\xi) =  \ol{g}_i(\xi)$ for all $\xi \in \Xi$. 
}

\begin{table*}[t]
\setlength{\tabcolsep}{5.4pt}
\renewcommand{\arraystretch}{1.2} 
\blueF{\begin{tabular*}{7.15in}{lllcccccccccccccccccccc} 
\toprule 
& & & &   \multicolumn{2}{c}{Inflexible Gen.}  & & \multicolumn{5}{c}{Intermittent Gen.} & &  \multicolumn{3}{c}{Flexible Gen.} & &\multicolumn{4}{c}{Load Buses}     \\  
\multicolumn{1}{c}{Parameters} & & Units & & 1 & 7 & & 2 &  3&  4&  5& 6&   & 8 & 9 & 10 & & 11 & 12 & 13 & 14  \\
\cmidrule{1-1} \cmidrule{3-3} \cmidrule{5-6} \cmidrule{8-12} \cmidrule{14-16}  \cmidrule{18-21}
$\alpha_i $ & & \$/\text{MW}  &  & 30 & 30 & & 0 & 0 & 0 & 0 & 0    &  &  50 & 50 & 50 & & -- & --  & -- & --  \\          
$\real\{d_i\}$ & & MW & &  0 & 21.7 & & 87.8 & 0 & 9.0 &  3.5 & 13.5  & & 94.2 & 11.2 & 0 & & 7.6 & 29.5 & 6.1 & 14.9   \\
$\imag\{d_i\}$ & & MVAR & & 0 & 12.7 & & -3.9  & 0 & 5.8 & 1.8 & 5.8  & &  19 &  7.5 & 0  & & 1.6 & 16.6 & 1.6 & 5 \\
$\real\{g^{\max}_i\}$ & & MW & &  100 & 100 & & 30  & 30  & 30 &30  & 30   & & 80  &  50 &  40 & & -- & --  & -- & --  \\
$\imag\{g^{\max}_i\}$ & & MVAR & & 10  & 50 & & $\dagger$ & $\dagger$  & $\dagger$   &$\dagger$   & $\dagger$     & & 40 & 24 & 24 & & -- & --  & -- & --   \\
$\real\{g^{\min}_i\}$& & MW & & 0  & 0 & & 0 & 0  & 0   & 0  & 0    & & 0 &0 & 0 & & -- & --  & -- & --   \\
$\imag\{g^{\min}_i\}$& & MVAR & &  0 & -40 & &  $\dagger$ & $\dagger$  & $\dagger$   &$\dagger$   & $\dagger$    & &  0 & -6 & -6 & & -- & --  & -- & --   \\
$ \real\{r_i^{\max}\}=-\real\{r_i^{\min}\}$  & & MW & &0 & 0 & & $\infty$  & $\infty$ & $\infty$ & $\infty$ &$\infty$  & & $\infty$ &$\infty$ & $\infty$ & & -- & --  & --  & -- \\               
$ \imag\{r_i^{\max}\}=-\imag\{r_i^{\min}\}$  & & MVAR & &0 & 0 & & $\infty$  & $\infty$ & $\infty$ & $\infty$ &$\infty$  & & $\infty$ &$\infty$ & $\infty$ & & -- & --  & -- & -- \\ 
Bus index (per \cite{testcases})   & & & & 1  & 2 & & 4 & 7& 10 &11&13 & &3&6 & 8 & & 5 & 9 & 12 & 14  \\    \bottomrule           
\end{tabular*}}
\caption{ \blueF{(IEEE 14-bus system). \  Specification of each generator's location, marginal cost, and constraint parameters. The $\dagger$ symbol indicates that the corresponding value in the table is determined by equation \eqref{eq:q_max}.} 
} 
\label{table:2a}
\vspace{-0.1cm}
\end{table*}

\section{Sequential Convex Approximation Method} \label{sec:mainresults_III}  
In what follows, we describe a recursive method that builds upon our previous development to generate a sequence of cost-improving convex inner approximations to the RAC-OPF problem. Let $v_0 \in \Cn$ be a feasible solution  to the \blueF{zero-recourse} AC-OPF problem \eqref{prog:SDR_OPF_polar1},  and define the matrix $V_0 \coloneqq v_0e_1^T$. 
Consider a recursion of the form
\begin{align} \label{eq:iter_alg}
	(x_{t+1},V_{t+1},\Lambda_{t+1}) \   \in    \underset{(x,V,\Lambda) \in \Feas{V_t}}{\argmin} \   \tr(MH_0(V,V_t)).
\end{align}
Here, $\Feas{V_t}$ denotes the feasible set of problem $\Pcal_{\rm III}(V_t)$, which is parameterized by the matrix $V_t$. The recursive algorithm \eqref{eq:iter_alg} can be interpreted as a successive convex majorization-minimization method. \blueF{An optimal solution $(x_t, V_t)$ associated with each step $t$ of the recursive method can be mapped to a feasible solution $(x_t, v_t(\bm{\xi}))$ of the RAC-OPF problem \eqref{prog:SDR_OPF_polar}, where $v_t(\bm{\xi}) = V_t \bm{\xi}$ is an affine recourse policy in $\bm{\xi}$.}
Proposition \ref{prop:Algorithm_properties} establishes two important properties of the recursive method. The recursive method is  (i)  guaranteed to yield a nonempty convex inner approximation to the RAC-OPF problem at each step in the recursion, and is (ii)  guaranteed to generate a sequence of \emph{feasible} dispatch policies for the RAC-OPF problem with nonincreasing costs. 
The proof  of Proposition \ref{prop:Algorithm_properties} can be found in Appendix B-C.
 \begin{propositio}\label{prop:Algorithm_properties} Let $v_0 \in \Cn$ be a feasible solution  to the \blueF{zero-recourse} AC-OPF problem \eqref{prog:SDR_OPF_polar1},  and define the matrix $V_0 \coloneqq v_0e_1^T$.  Let $\{x_{t},V_{t},\Lambda_{t}\}_{t=1}^{\infty}$ denote the sequence of solutions generated by the recursion in \eqref{eq:iter_alg}. The following properties hold for each step $t$ of the recursion.
 \vspace{.1in}
	\begin{enumerate}[(i)]\itemsep0.5em
		\item Nonemptiness: $\Feas{V_t} \neq \emptyset$. \label{prop:res_1}
		\item  Cost montonicity:  $\tr(MV_{t}^\herm A_0V_{t}) \le   \tr(MV_{t-1}^\herm A_0V_{t-1})$.  \label{prop:res_2}
	\end{enumerate}
\end{propositio}
 \blueF{We also remark that it is also possible to establish convergence of the recursive method \eqref{eq:iter_alg} to a stationary point of the nonconvex problem \ref{prog:two_stage_Sprocedure1} using existing techniques from the literature  \cite{lanckriet2009convergence, lipp2016variations}.}

\section{Case Study} \label{sec_numerical_studies}
We now  illustrate the effectiveness of the proposed optimization method on the \blueF{WSCC 9-bus \cite{sauer1998power} and IEEE 14-bus \cite{testcases} test systems} with different levels of renewable resource penetration and uncertainty.  We refer the reader to the aforementioned references for the complete specification and single-line diagrams of the test systems considered. All modifications made to the original WSCC 9-bus and IEEE 14-bus test systems are summarized in Tables \ref{table:1a} and \ref{table:2a}, respectively. \blueF{Additionally, for each test system considered, we permit load shedding at non-generator buses (i.e., load buses), where a reduction in load relative to the specified level is penalized according to the value of lost load (VOLL).\footnote{We refer the reader to Appendix \ref{app:Load_Shedding} for a complete specification of the RAC-OPF problem with load-shedding.} We set the VOLL equal to \$4,000/MWh.}

\subsection{Renewable Generator Model} \label{sec:rengenmodel}
The real-time generating capacity of renewable generators $i \in \{2, \dots, 6\}$ represents the only source of uncertainty in the power system being considered. Accordingly, we set the dimension of the random vector to $k=6$, and let the $i^{{\rm th}}$ element of the random vector $\bm{\xi}$ represent the \emph{maximum active power} available to generator $i$ in real-time. In other words, 
\begin{align} \label{eq:num_con}
\real\{\underline{g}_i(\bm{\xi})\} =  0  \ \  \text{and} \ \  \real\{\overline{g}_i(\bm{\xi})\} = \bm{\xi}_i,
\end{align}
for $ i = 2, \dots, 6$.
It will be convenient to our numerical analyses in the sequel to express the random vector $\bm{\xi}$ as an affine function of a zero-mean random vector $\bm{\delta}$ that is uniformly distributed over a unit ball. We define this relationship  according to 
\begin{align*}
\bm{\xi} \coloneqq \mu + \sigma \bm{\delta},
\end{align*}
where the random vector $\bm{\delta}$ is assumed to have support 
\begin{align*}
\Delta  \coloneqq \{\delta \in \R^k \ | \ \delta_1 = 0, \ \|\delta\|_2\le 1 \}.
\end{align*}
It follows that the random vector $\bm{\xi}$ has support given by
\begin{align*}
\Xi =\left\{ \xi \in \R^k  \  \left| \    \xi-\mu \in \sigma\Delta \right. \right\}.
\end{align*} 
Here, $\mu \in \R^k$ and $\sigma \in \R_+$  represent \emph{location} and \emph{scale parameters}, respectively.  In the following study, we set 
$\mu_i = 15$ MW for each renewable generator $i \in \{2, \dots, 6\}$. Qualitatively, the larger the scale parameter $\sigma$, the larger the a priori uncertainty in the real-time generating capacity of the renewable generators. 
The location and scale parameters are chosen in such a manner as to ensure that $\bm{\xi}_i$ respects the nameplate active power capacity limits   for each renewable generator $i$  specified in Table \ref{table:1a}. We also require that $\mu_1 = 1$ to maintain consistency with our original uncertainty model in Section \ref{sec:II_Uncertainty}.  Finally, under the assumption that $\bm{\delta}$ has a uniform distribution, it is straightforward to show that the random vector $\bm{\xi}$ has a second-order moment matrix given by
\begin{align*}
M  = \mu \mu^T + \left(\frac{\sigma^2}{k+1} \right)
\left[
\begin{array}{c|c}
0 &  \\ \hline
 & I_{k-1}
\end{array}\right].
\end{align*}

\begin{figure*}[t] 
\centering
\begin{subfigure}{0.32\textwidth}
\begin{overpic}[scale=0.32]{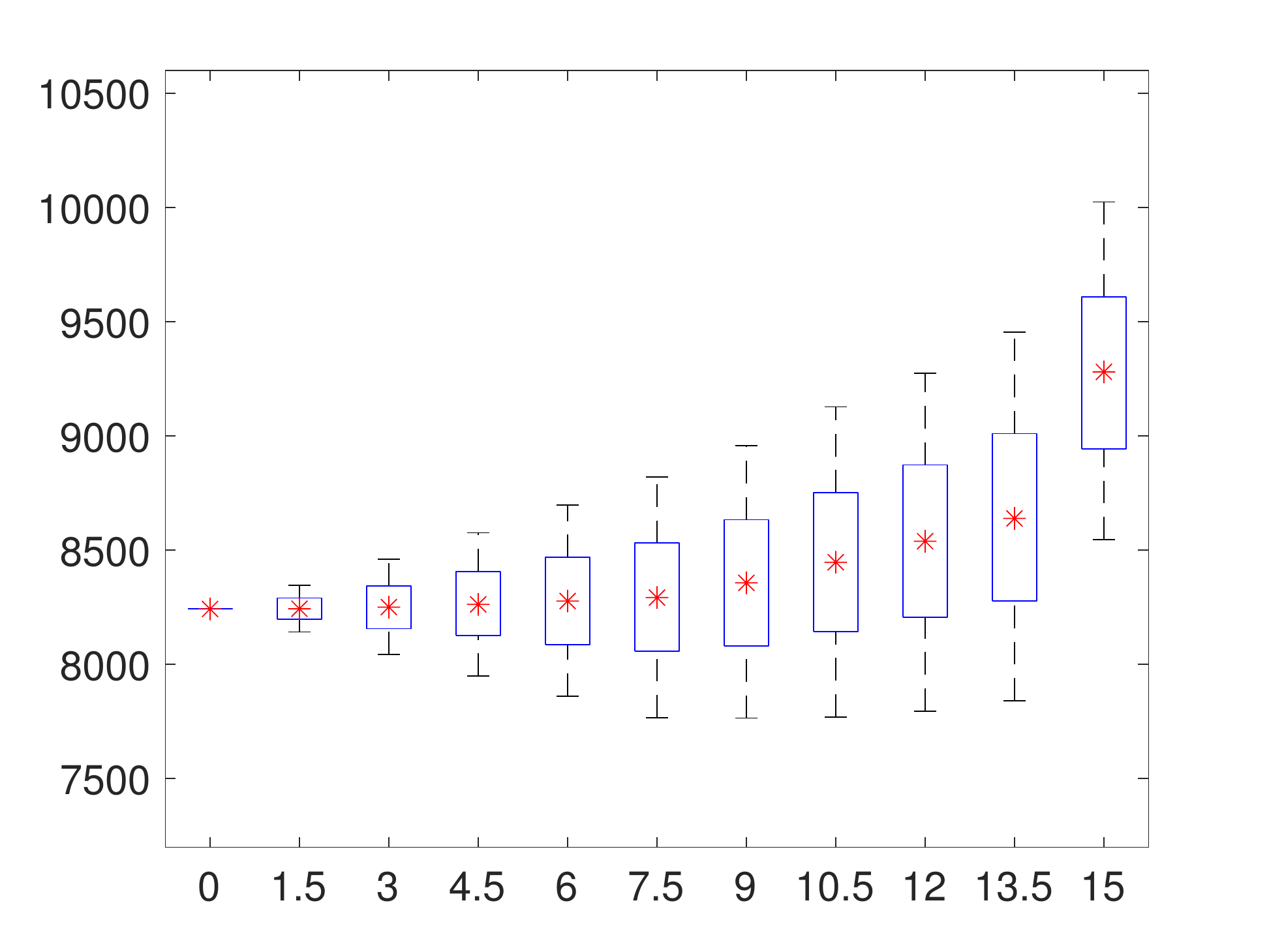}
\put (47,-1.0) { \footnotesize{$\sigma$}}
\put (22,-6.9) { \footnotesize{(a) Cost Distribution (9-bus)}}
\put (-2,40.5) { \footnotesize{\$} }
\end{overpic}
\label{fig:subfig_c}
\end{subfigure}
\begin{subfigure}{0.32\textwidth}
\begin{overpic}[scale=0.32]{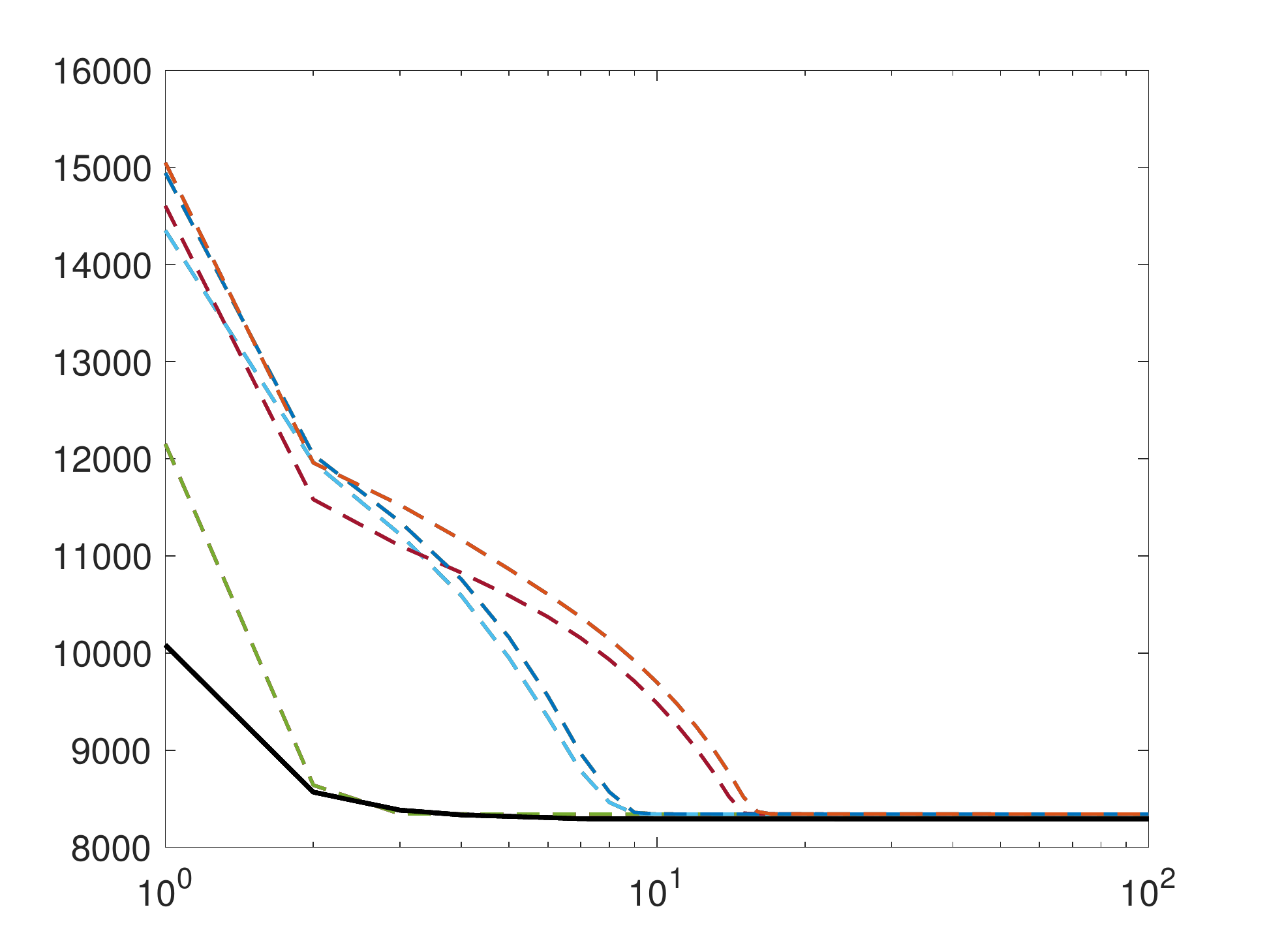}
\put (45,-1.5) { \footnotesize{$t+1$} }
\put (20,-6.9) { \footnotesize{(c) Recursive Method (9-bus)}}
\put (0,40.5) { \footnotesize{\$} }
\end{overpic}
\label{fig:subfig_a}
\end{subfigure}
\begin{subfigure}{0.32\textwidth}
\begin{overpic}[scale=0.32]{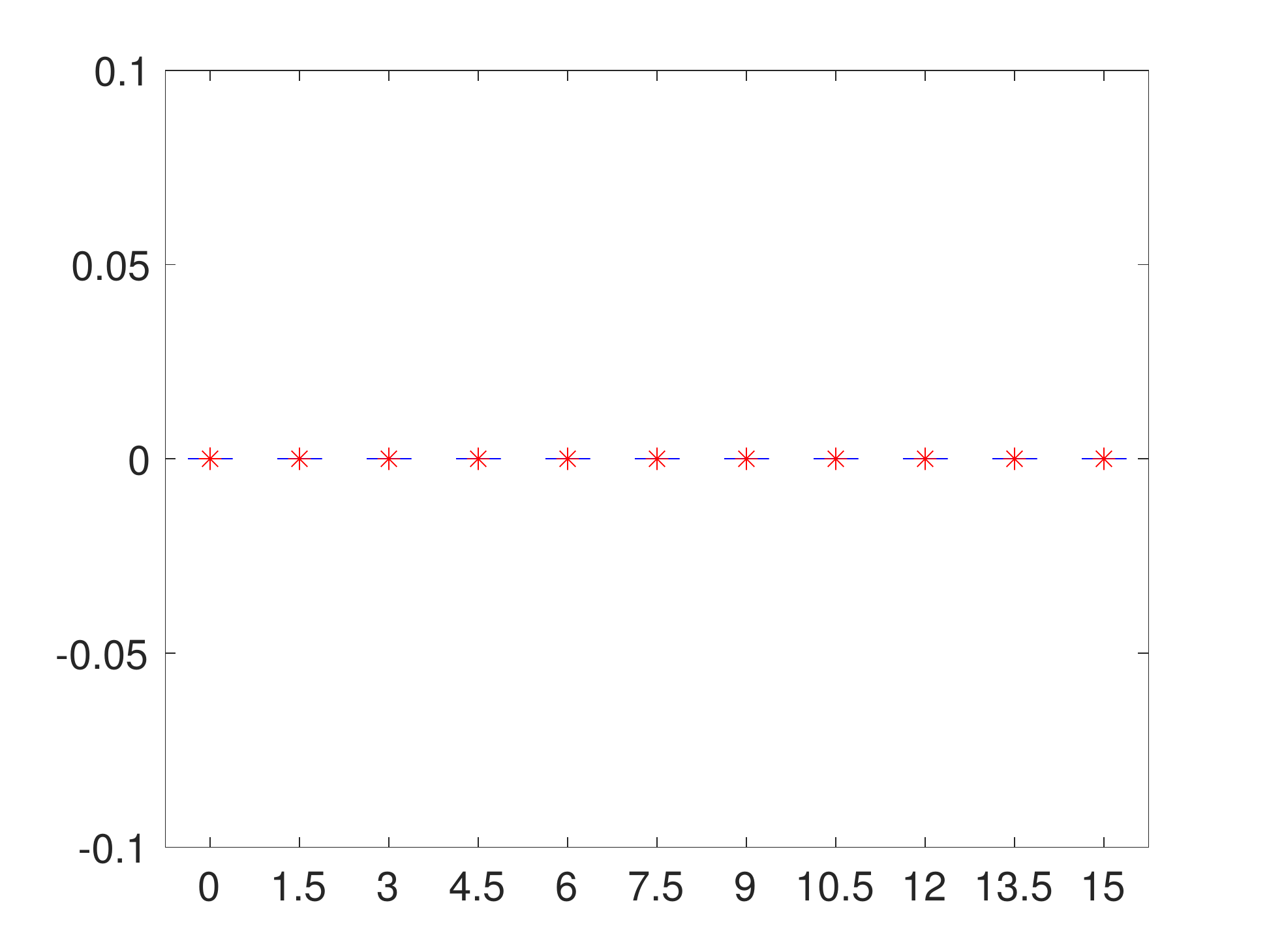}
\put (47,-1.0) { \footnotesize{$\sigma$}}
\put (17,-6.9) { \footnotesize{(e) Load Shed Distribution (9-bus)}}
\put (-1,40.5) { \footnotesize{MW} }
\end{overpic}
\label{fig:subfig_e}
\end{subfigure}\\[0.3em]
\begin{subfigure}{0.32\textwidth}
\begin{overpic}[scale=0.32]{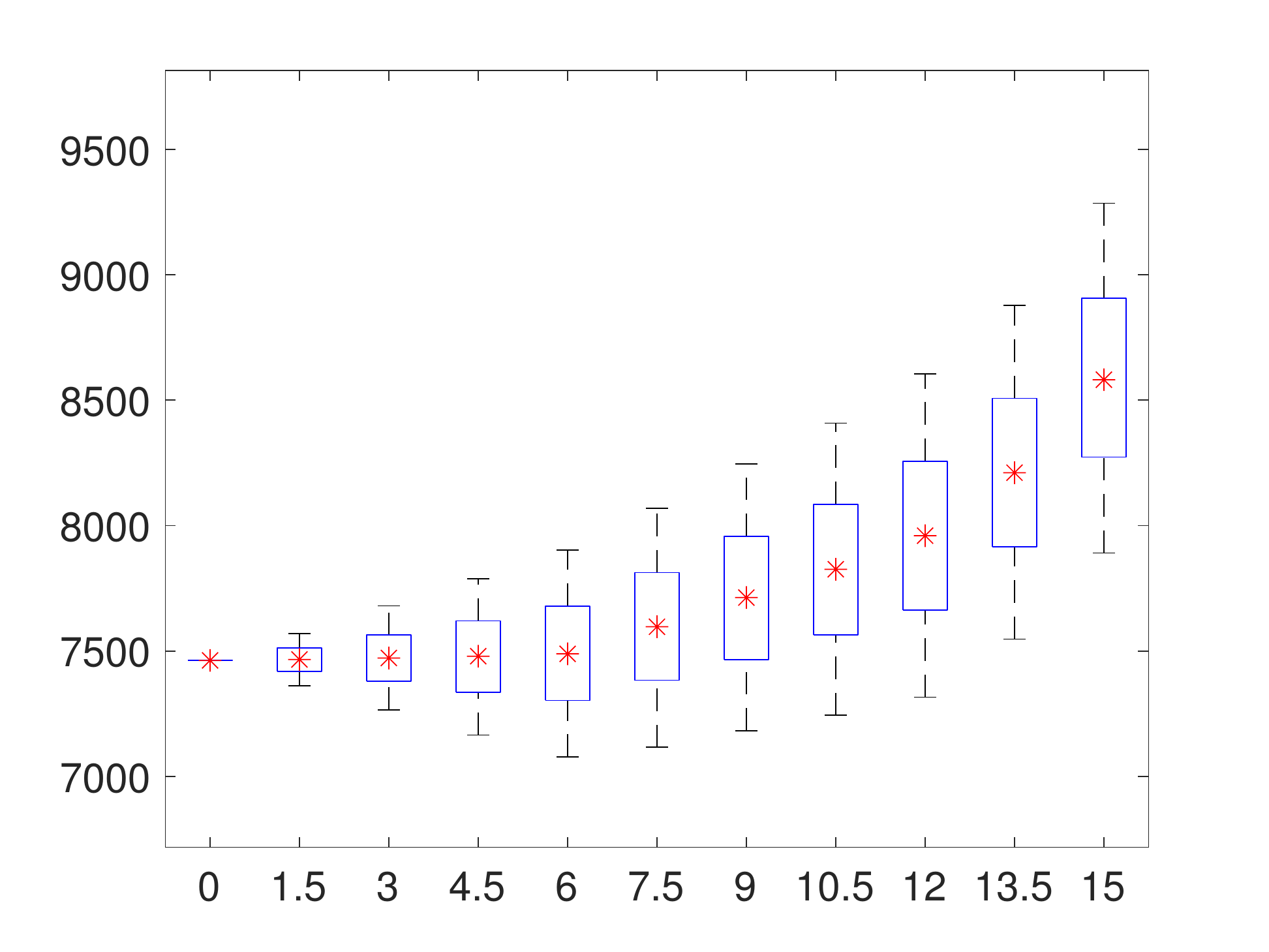}
\put (47,-1.0) { \footnotesize{$\sigma$}}
\put (22,-6.9) { \footnotesize{(b) Cost Distribution (14-bus)}}
\put (-2,40.5) { \footnotesize{\$} }
\end{overpic}
\label{fig:subfig_d}
\end{subfigure} 
\begin{subfigure}{0.32\textwidth}
\begin{overpic}[scale=0.32]{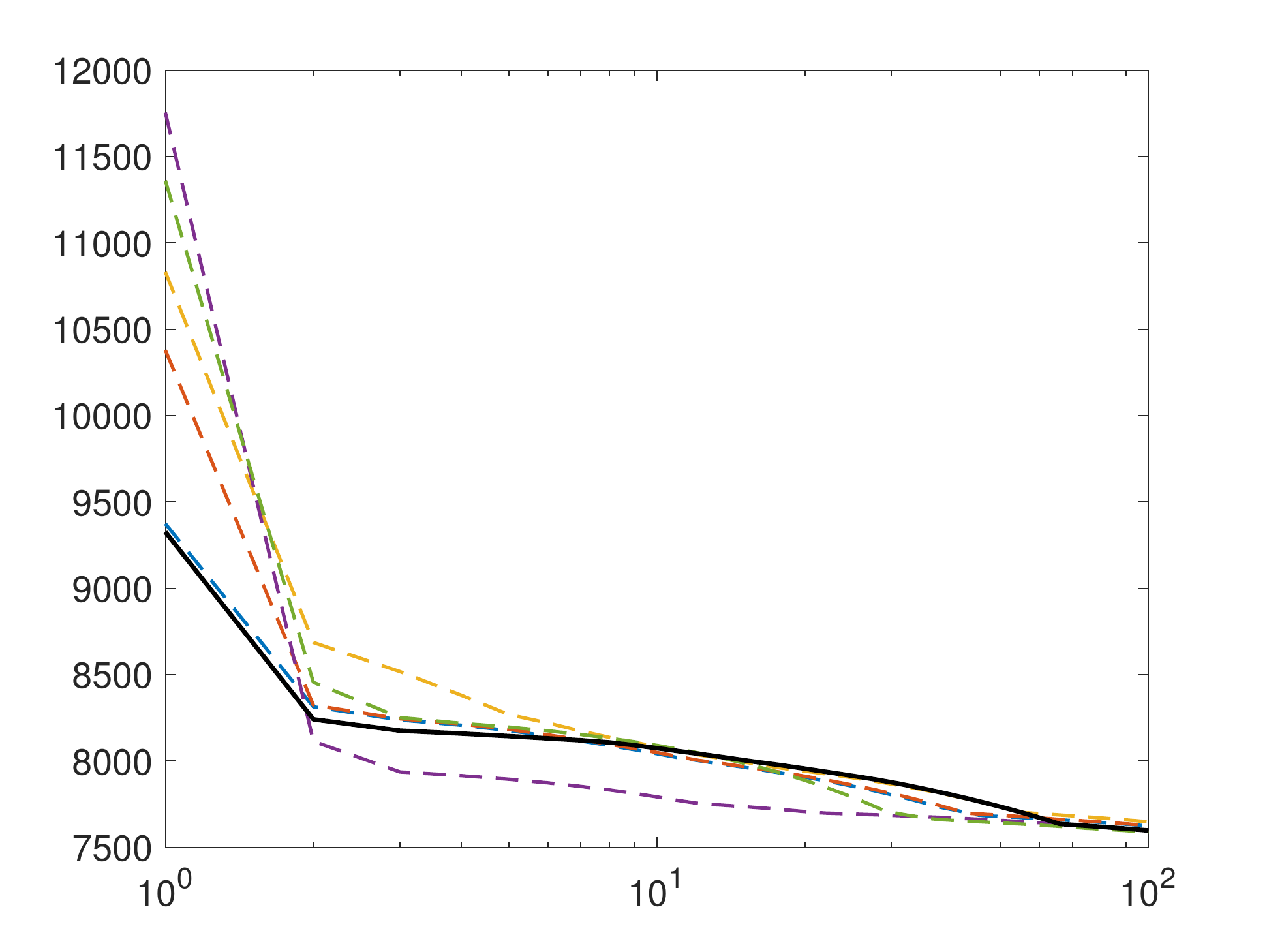}
\put (45,-1.5) { \footnotesize{$t+1$} }
\put (20,-6.9) { \footnotesize{(d) Recursive Method (14-bus)}}
\put (-2,40.5) { \footnotesize{\$} }
\end{overpic}
\label{fig:subfig_b}
\end{subfigure} \ \ 
\begin{subfigure}{0.32\textwidth}
\begin{overpic}[scale=0.32]{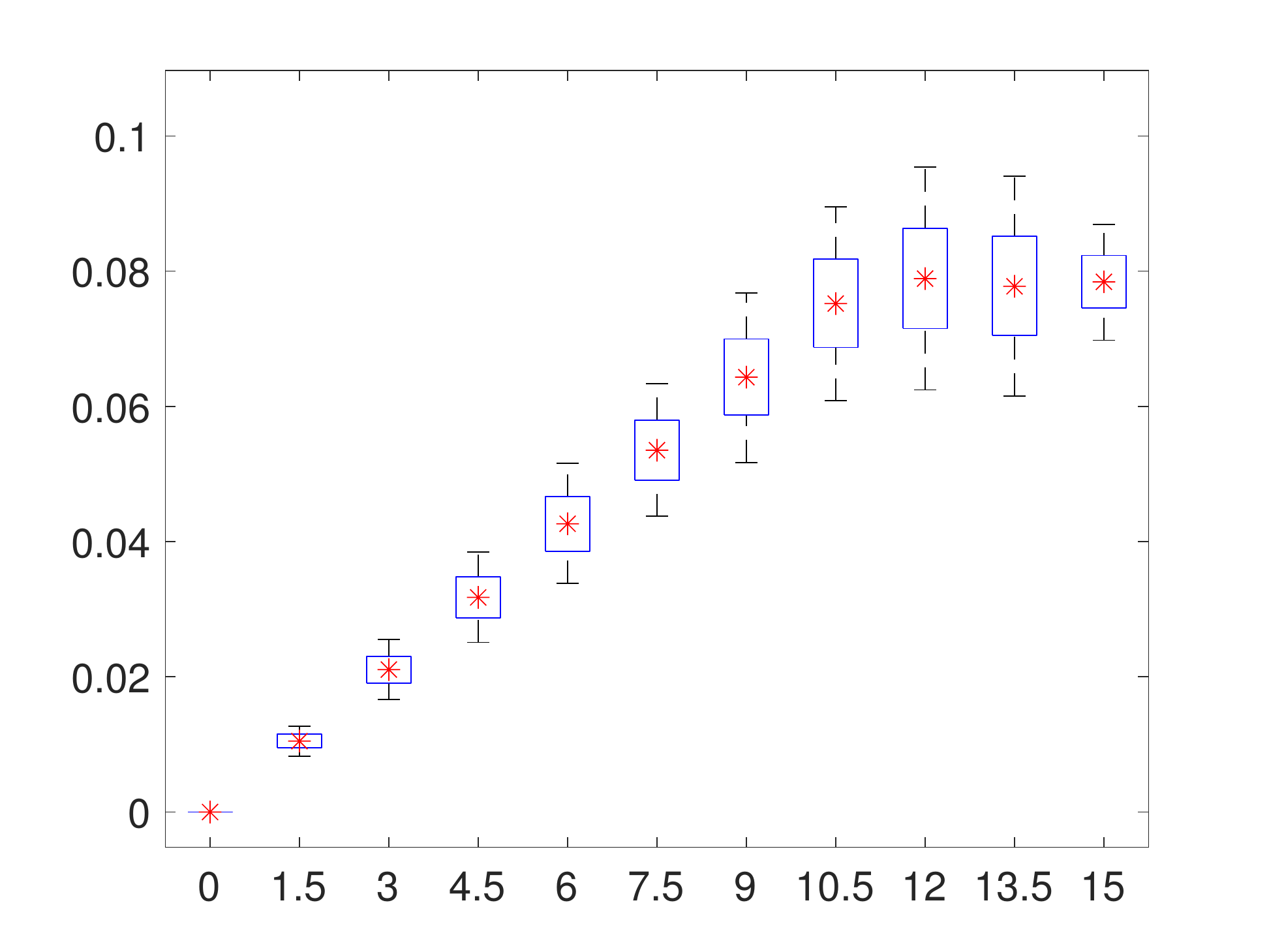}
\put (47,-1.0) { \footnotesize{$\sigma$}}
\put (17,-6.9) { \footnotesize{(f) Load Shed Distribution (14-bus)}}
\put (-5.8,40.5) { \footnotesize{MW} }
\end{overpic}
\label{fig:subfig_f}
\end{subfigure}

\caption{
The top (bottom)  row  of figures is associated with the WSCC 9-bus (IEEE 14-bus) test system. 
\textbf{Figures (a)-(b)} depict the expected generation cost (red star) and empirical confidence intervals incurred by the affine dispatch policy returned by the recursive algorithm \eqref{eq:iter_alg} versus the scale parameter $\sigma$.   The empirical confidence intervals are estimated using 10,000 independent realizations of the underlying random vector $\bm{\xi}$.  The box depicts the  interquartile  range, while the lower and upper whiskers extend to the 5\% and 95\% quantiles, respectively. 
 \textbf{Figures (c)-(d)} plot (as colored  dashed lines) the expected generation cost (for $\sigma = 7.5$)  incurred at  each step of the recursive algorithm \eqref{eq:iter_alg}  for five randomly generated initial conditions $V_0$. The solid black line represents the expected cost trajectory returned by the recursive algorithm given the initial condition $V_0$  generated by the zero-recourse AC-OPF problem \eqref{prog:SDR_OPF_polar1}.
 \textbf{Figures (e)-(f)}  depict the  empirical confidence intervals for the total active power load that is shed at non-generator buses under the affine dispatch policy returned by algorithm \eqref{eq:iter_alg} versus the scale parameter $\sigma$. The empirical confidence intervals are generated using 10,000 independent realizations of the random vector $\bm{\xi}$ for each value of $\sigma$. The box depicts the  interquartile  range, while the lower and upper whiskers of the confidence intervals extend to the minimum and maximum values of the total active power load shed.
}
\label{fig:iterations}
\end{figure*}
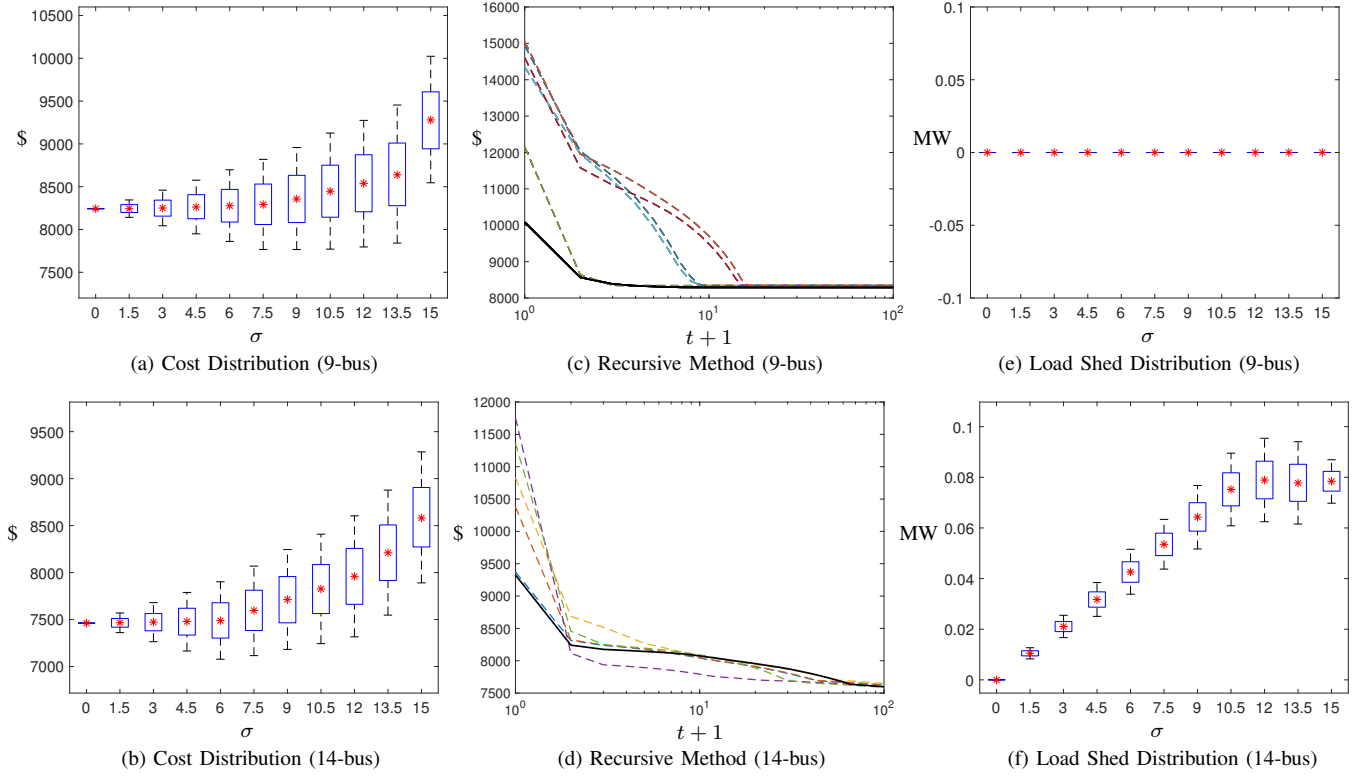

Renewable energy resources, like wind and solar, employ power electronic inverters, which can produce and absorb reactive power. The limits on the maximum and minimum amount of reactive power that can be injected by  a renewable generator are determined by its inverter's apparent power capacity, which we denote by $s_i^{\max} \in \R_+$  for each renewable generator $i$. It follows that the real-time complex power injection of each renewable generator $i$ must  satisfy a capacity constraint of the form 
\begin{align}
|g_i(\bm{\xi})|  \le s_i^{\max}.  \label{eq:apparent}
\end{align}
As the slight oversizing of a renewable generator's apparent power rating is standard in practice, we set  $s_i^{\max} = 1.05\real\{g_i^{\max}\}$  for each renewable generator $i$.
In order to ensure that Assumption \ref{ass:lin} is satisfied, we enforce a more conservative form of the real-time apparent power capacity constraint \eqref{eq:apparent} by setting the real-time reactive power limits for each renewable generator $i$ according to
\begin{equation}
\begin{aligned}\label{eq:q_max}
 	\imag\{\overline{g}_i(\bm{\xi})\} & = -\imag\{\underline{g}_i(\bm{\xi})\} \\ &= \inf_{\xi \in \Xi} \ \sqrt{(s_i^{\max})^2 - \xi_{i}^2} \\
 	&= \sqrt{(s_i^{\max})^2 - (\mu_i+\sigma)^2} .
\end{aligned}
\end{equation}
\blueN{For simplicity, we also fix the nameplate reactive power limits of each renewable generator according to 
$\imag\{g^{\max}_i\} =\imag\{\overline{g}_i(\bm{\xi})\}$ and $\imag\{g^{\min}_i\} = \imag\{\underline{g}_i(\bm{\xi})\}.$}
The reactive power limits specified in \eqref{eq:q_max} specify the range of reactive power injections that are guaranteed to be available to a renewable generator in real-time, regardless of the active power supplied. 
Using the real-time active and reactive power capacity limits specified in \eqref{eq:num_con} and \eqref{eq:q_max}, respectively, it is straightforward to construct matrices $\ul{G}, \ol{G} \in \C^{n \times k}$ such that Assumption \ref{ass:lin} is satisfied.

\blueN{Finally, it is worth noting that the conservative generation  limits $(\gamma_i^{\min}, \gamma_i^{\max})$ defined in \eqref{eq:gammata} admit closed-form expressions for the system being considered. They are given by
\begin{equation}
\begin{aligned}
\gamma_i^{\min}&=  \hspace{.2in} 0   \hspace{.3in}  - \ &\mathbf{j}\sqrt{(s_i^{\max})^2 - (\mu_i+\sigma)^2}, \\
\gamma_i^{\max} &= (\mu_i-\sigma)   \ \ +   &\mathbf{j}	 \sqrt{(s_i^{\max})^2 - (\mu_i+\sigma)^2}.
\end{aligned}
\end{equation}}

\vspace{-0.2cm}
\subsection{Numerical Analyses and Discussion}
We begin by examining the sensitivity of the generation cost incurred under the affine recourse policies that we propose to uncertainty in renewable supply. We do so by varying the scale parameter $\sigma$ from 0 to  15 in increments of 1.5,
 while keeping all other problem parameters fixed. It is worth noting  that for $\sigma= 0$, there is  no a priori uncertainty in the renewable supply, and the RAC-OPF problem  \eqref{prog:SDR_OPF_polar}  reduces to the zero-recourse AC-OPF problem \eqref{prog:SDR_OPF_polar1}. \blueN{In addition, for $\sigma = 0$, we are able to verify that the solutions we compute are in fact an optimal solution by using the semidefinite relaxation of the AC-OPF problem described in \cite{low2014convexII}. }
For each value of $\sigma$ that we consider, we calculate an affine recourse policy according to the recursive algorithm specified in Eq.  \eqref{eq:iter_alg}.\footnote{\blueN{All numerical analyses were carried
out in Matlab and semidefinite programs were solved using SDPT3 \cite{toh1999sdpt3}. The machine used to
to solve the problems has a 3.1GHz  Intel dual-core with 16GB of RAM.}} We initialize the recursion with a feasible solution to the zero-recourse AC-OPF problem  \eqref{prog:SDR_OPF_polar1}, which  we compute using the Matpower interior point solver \cite{zimmerman2011matpower}. \blueN{The recursive algorithm terminates when either one of the two following conditions hold: (i) the difference
between the optimal values of successive iterations is less than $10^{-4}$, or (ii) the total number of iterations exceeds 500.}

We plot the expected generation cost (and empirical confidence intervals) incurred by the affine dispatch policy returned by algorithm \eqref{eq:iter_alg} versus the scale parameter $\sigma$ for the 9-bus and 14-bus systems in Figures \ref{fig:iterations}(a) and \ref{fig:iterations}(b), respectively. First, notice that, for each test system, the expected generation cost increases monotonically with the scale parameter. Such behavior is to be expected,  as  larger values of $\sigma$ correspond to larger uncertainty sets $\Xi$. It is also worth noting the `spread' in the cost distribution induced by the dispatch policies that we compute also increases with $\sigma$. That is to say, renewable energy resources with a large variance in their real-time generating  capacity will result in a larger variance in total generating costs. Such behavior is a consequence of the risk neutrality inherent to the expected cost criterion that we treat in this paper. 

\blueF{
We examine the convergence behavior of the recursive algorithm  \eqref{eq:iter_alg} for the  9-bus and 14-bus test systems in Figures \ref{fig:iterations}(c) and \ref{fig:iterations}(d), respectively. For each system, we set $\sigma = 7.5$ and plot (as dashed lines) the expected generation cost  incurred at  each step of the recursive algorithm \eqref{eq:iter_alg} for five randomly generated initial conditions $V_0$.\footnote{We initialize the recursive algorithm \eqref{eq:iter_alg} with randomly generated matrices $V_0$ by solving the zero-recourse AC-OPF problem \eqref{prog:SDR_OPF_polar1} with with randomly generated marginal-cost parameters $\{\alpha_i\}$. Specifically, the generators' marginal-cost parameters $\{\alpha_i\}$ are sampled independently and uniformly at random from the interval $[0,50]$.}  In both  Figures \ref{fig:iterations}(c) and \ref{fig:iterations}(d), the solid black lines represent the expected cost trajectory returned by the recursive algorithm given the initial condition $V_0$  generated by the zero-recourse AC-OPF problem \eqref{prog:SDR_OPF_polar1} under the \emph{true} marginal-cost parameters specified in Tables \ref{table:1a} and \ref{table:2a}. It is interesting to note that the recursive algorithm converges to the same optimal value regardless of the initial condition. 
 Additionally, the numerical results in Figures \ref{fig:iterations}(c)-(d) are consistent with the guarantees of Proposition \ref{prop:Algorithm_properties}, which ensures that the recursive algorithm \eqref{eq:iter_alg} will yield a sequence of feasible dispatch policies with nonincreasing costs.  We also note that each iteration for the 9-bus (14-bus)  system took 11.67 seconds (130.4 seconds) to complete on average.}

\blueF{ In Figures \ref{fig:iterations}(e) and \ref{fig:iterations}(f), we plot empirical confidence intervals for the total active power load that is shed at non-generator buses under the affine dispatch policy returned by algorithm \eqref{eq:iter_alg} versus the scale parameter $\sigma$ for the 9-bus and 14-bus systems, respectively.    The confidence intervals are generated using 10,000 independent realizations of the random vector $\bm{\xi}$ for each value of $\sigma$. The lower and upper whiskers of the confidence intervals extend to the minimum and maximum values of the total active power load shed. For the 9-bus system, the affine dispatch policy never sheds active power load, as depicted in Figure \ref{fig:iterations}(e). That is to say, the allowance of load shedding does not result in any reduction of load under the affine dispatch policies that we compute. 
For the 14-bus system, we empirically observe that only a small percentage of total active power load is shed under the affine dispatch policies that we compute.  In particular, the total active power load shed never exceeds 0.16\% of the total active power load  at non-generator buses.  Lastly, we remark that for both the 9-bus and 14-bus systems, we empirically observe that the affine dispatch policy never sheds reactive power load.}

\section{Conclusion} \label{sec:conc} 
In this paper, we formulate the robust AC optimal power flow (RAC-OPF) problem as a two-stage robust optimization problem with recourse. \blueF{Under the assumption that there is adjustable generation or load at every bus in the power transmission network,}  we provide a technique to construct a convex inner approximation of RAC-OPF in the form of a semidefinite program. 
In particular, the inner approximation is obtained by: (i) restricting the set of admissible recourse policies to be affine in the uncertain variables, (ii) approximating the semi-infinite constraint set by a sufficient set of finitely-many constraints, and (iii) approximating the indefinite quadratic constraints by majorizing convex quadratic constraints.
Its solution yields an affine recourse  policy that is \emph{guaranteed to be feasible} for RAC-OPF.  
 In addition, we provide an iterative optimization algorithm that generates a sequence of feasible affine recourse policies with nonincreasing costs.
 
 There are several interesting directions for future research. 
 First, affine recourse policies are likely to be suboptimal for the RAC-OPF problem. Thus, it would be interesting to investigate the design of convex relaxations for RAC-OPF to enable the tractable calculation of  lower bounds on the optimal value of RAC-OPF. 
 Such lower bounds can, in turn, be used to bound the suboptimality incurred by the feasible affine policies proposed in this paper. \blueF{Second, the convex inner approximation technique developed in this paper relies explicitly on the rather limiting assumption that there is adjustable generation or load at every bus in the power transmission network. In order to accommodate the treatment of more general power systems, it will be important to relax this assumption, while preserving the robust feasibility guarantees developed in this paper.
  It would also be of interest to extend the techniques developed in this paper  to accommodate discrete decision variables (e.g., unit commitment decisions) in the RAC-OPF problem.}

\vspace{-0.3cm}

\bibliographystyle{IEEEtran}	
\bibliography{references_bib}

\begin{appendices}
\label{sec:appendix}
\vspace{-0.1cm}
\section{Concise Reformulation of RAC-OPF} \label{app:matricesAi}

Define matrices $\Phi_i, \Psi_i \in \Hn$, for all $i \in \Vcal$, and a matrix $E \in \R^{4n \times 2n}$ as follows:
\begin{align*}
	 \Phi_i \coloneqq \frac{ S_i + S_i^\herm }{ 2},  \ \   \Psi_i \coloneqq \frac{ S_i - S_i^\herm }{\mathbf{j}2}, \ \ E \coloneqq \bmat{I_{2n} & -I_{2n}}^T.
\end{align*}
In addition, let $f \in \R^{4n}$ be a vector given by
\begin{align*}
f \coloneqq \left[\real\{g^{\max}\}^T \ \ \imag\{g^{\max}\}^T \ -\real\{g^{\min}\}^T \   -\imag\{g^{\min}\}^T \right]^T
\end{align*}

The RAC-OPF problem \eqref{prog:SDR_OPF_polar} can be reformulated as follows: 
\begin{alignat}{4}\label{prog:SDR_OPF_pola11}
		& \underset{}{\text{minimize}}  \ 
		 \Exp \left[ \sum_{i=1}^n \alpha_i  v(\bm{\xi})^\herm \Phi_i  v(\bm{\xi}) \right] &&\hspace{-2.1cm}+\sum_{i=1}^n\alpha_i\real\{\dnom_i\}  \nonumber \\[0.6em]
		& \text{subject to}  \\
		&	   x \in \R^{2n}, \  \ v \in \Lb^2_{k,n}   \nonumber \\[0.5em]
		&     v(\xi)^\herm \Phi_i v(\xi)  \le  \real\{ e_i^T\ol{G}-\dnom_ie_1^T\}\xi, && i \in \Vcal \nonumber \\[0.5em]
		&    v(\xi)^\herm (-\Phi_i) v(\xi)  \le  \real\{\dnom_ie_1^T- e_i^T\ul{G}\}\xi, && i \in \Vcal \nonumber  \makebox[30pt][r]{\smash{\raisebox{-6.3\baselineskip}{$\left.\rule{0pt}{8.3\baselineskip}\right\} \rotatebox{90}{$\hspace{-0.5cm}\forall \ \xi \in \Xi.$} \quad $}}} \\[0.5em]
		&    v(\xi)^\herm \Psi_i v(\xi)  \le  \imag\{ e_i^T\ol{G}-\dnom_ie_1^T\}\xi, && i \in \Vcal  \nonumber  \\[0.5em]
		&  v(\xi)^\herm (-\Psi_i) v(\xi)  \le \imag\{\dnom_ie_1^T- e_i^T\ul{G}\}\xi, && i \in \Vcal  \nonumber  \\[0.5em]
		&    v(\xi)^\herm \Phi_i v(\xi) - e_i^Tx  \le  \real\{r_i^{\max}-\dnom_i\}e_1^T\xi, && i \in \Vcal   \nonumber \\[0.5em]
		&    v(\xi)^\herm (-\Phi_i) v(\xi) + e_i^Tx  \le  \real\{\dnom_i-r_i^{\min}\}e_1^T\xi, && i \in \Vcal  \nonumber \\[0.5em]
		&    v(\xi)^\herm \Psi_i v(\xi) - e_{n+i}^Tx  \le  \imag\{r_i^{\max}-\dnom_i\}e_1^T\xi, && i \in \Vcal \nonumber \\[0.5em]
		&    v(\xi)^\herm (-\Psi_i) v(\xi) + e_{n+i}^Tx  \le  \imag\{\dnom_i-r_i^{\min}\}e_1^T\xi, \  && i \in \Vcal  \nonumber \\[0.5em]
		&   v(\xi)^\herm e_ie_i^Tv(\xi)  \le  (v_i^{\max})^2e_1^T\xi, && i \in \Vcal \nonumber \\[0.5em]
		&   v(\xi)^\herm (-e_ie_i^T)v(\xi)  \le  -(v_i^{\min})^2e_1^T\xi,  && i \in \Vcal \nonumber \\[0.5em]		
		&    v(\xi)^\herm P_{ij} v(\xi)  \le   \ell_{ij}^{\max}e_1^T\xi,  && \hspace{-0.64cm} (i,j)\in \Ecal \nonumber\\[0.5em]
		&    v(\xi)^\herm  (-P_{ij}) v(\xi)  \le   \ell_{ij}^{\max}e_1^T\xi,  && \hspace{-0.64cm} (i,j)\in \Ecal \nonumber\\[0.5em]
	   &   E x \le f. \nonumber 
\end{alignat}

\section{Proofs}
\subsection{Proof of Lemma \ref{lemma:properties_Hi}} \label{app:proof_lemma_prop_Hi}
\noindent \emph{Part \eqref{lemma:properties_Hi_a}:} \  Notice that 
\begin{align}  \label{eq:matrix_linear1}
V^\herm A_iV - H_i(V,Z) = (Z-V)^\herm A_i^-(Z-V) .
\end{align}
The desired result follows from the fact that $\Ans_i \preceq 0$.

\

\noindent \emph{Part \eqref{lemma:properties_Hi_b}:} \  Let $N \in \R^{k\times k}$ be a Cholesky factor of $M$, i.e., $M = NN^T$. The existence of $N$ is guaranteed as $M$ is assumed to be positive definite.
Since the matrix on the left-hand side of \eqref{eq:matrix_linear1} is Hermitian negative semidefinite, we must have
\[
	N^T ( V^\herm A_iV - H_i(V,Z))N   \preceq  0,
\]
And since the trace of any Hermitian negative semidefinite matrix is nonpositive, we obtain
\[
	\tr(N^T (V^\herm A_iV -H_i(V,Z))N)   \le  0.
\]
The linearity of the trace operator and its invariance under cyclic permutations gives the desired trace inequality. \qed

\subsection{Proof of Proposition \ref{prop:feasibility}} \label{app:proof_prop_feas}
Given a feasible solution $v_0$ to problem \eqref{prog:SDR_OPF_polar1}, we show that there exists $x_0 \in \R^{2n}$ and $\Lambda_0 \in \R^{m \times \ell}$ such that the point $(x_0,V_0,\Lambda_0)$ is feasible for problem \ref{prog:two_stage_Linearization}, where $V_0=v_0e_1^T$. Let $x_0 \in \R^{2n}$ be a vector given by
\[
	[x_0]_i \coloneqq \begin{cases}
			\real\{v_0^\herm S_i v_0 + d_i \}, & \text{if } 1 \le i \le n, \\[0.3em]
			\imag\{v_0^\herm S_{i-n} v_0 +d_{i-n} \}, & \text{if } n < i \le 2n,
			\end{cases}
\]
where $S_i \in \C^{n\times n}$ is defined in \eqref{eq:power_balance}. Since $v_0$ is a feasible solution to \eqref{prog:SDR_OPF_polar1}, it is easy to verify that $Ex_0 \le f$, where 
$E \in \R^{4n \times 2n}$ and $f \in \R^{4n}$ are defined in Appendix \ref{app:matricesAi}. Therefore, it suffices to show 
  that there exists $\Lambda_0 \le 0$ such that
\begin{align}\label{eq:proof_must_showLambdo0}
	H_i(V_0,V_0)-C_i+(b_i^T x_0)e_1e_1^T-\sum_{j=1}^\ell [\Lambda_0]_{ij}W_j \preceq 0,
\end{align}
for all $i \in [m]$. Since $V_0 = v_0e_1^T$, we have 
\begin{align} \label{eq:help_proof_feas}
	H_i(V_0,V_0) = V_0^\herm A_iV_0 = (v_0^\herm A_iv_0)e_1e_1^T,
\end{align}
for all $i \in [m].$
Therefore, \eqref{eq:proof_must_showLambdo0} takes the form 
\begin{align} \label{eq:help_proof_feas}
	(v_0^\herm A_iv_0+b_i^T x_0)e_1e_1^T -C_i - \sum_{j=1}^\ell [\Lambda_0]_{ij}W_j \preceq 0, 
\end{align}
for all $i \in [m]$. Fix an $i \in [m]$  and consider
the constraint
\[
	v^\herm A_iv + b_i^T x \le c_i^T \xi, \quad \forall \ \xi \in \Xi,
\]
in problem \ref{prog:two_stage_general}. This constraint can be equivalently
expressed as $v^\herm A_iv + b_i^T x \le \pi_i$, where $\pi_i$ is the optimal value of the following convex quadratically constrained quadratic program
\begin{equation}
\begin{alignedat}{4} \label{prog:equiv_myopic}
&   \underset{}{\text{minimize}} \quad
	& & \xi^T C_i\xi \nonumber \\
	& \text{subject to}  \quad   
	&&  \xi \in \R^{k} \\[0.2em]
	&&&  \xi^T W_j\xi \ge 0, \quad  \ j\in[\ell],\\
	&&& \xi_1=1.
\end{alignedat} 
\end{equation} 
In the above optimization problem, we have used the relation $\xi_1 = 1$ to rewrite $c_i^T \xi$ as $\xi^TC_i\xi$, where  $C_i \in \R^{k\times k}$ is defined in problem \ref{prog:two_stage_Sprocedure1}. 
By Asumption \ref{ass:slater}, there exists $\xi \in \Xi$ such that $\xi^T W_j\xi > 0$ for all $j \in [\ell]$.  This is a Slater condition, which guarantees strong duality to hold between the above optimization problem and its dual. Therefore,  $\pi_i = \tau_i$, where $\tau_i$ is the optimal value of its dual problem, which is given by
\begin{equation}
\begin{alignedat}{4} \label{prog:equiv_myopic}
&  \underset{}{\text{maximize}} \quad
	& & \rho_i + \gamma_i \nonumber \\
	& \text{subject to}      \quad
	&&  \rho_i \in \R, \ \gamma_i \in \R,  \ [\Lambda_0]_{ij} \in \R,  \ j \in [\ell] \\[0.5em]
	&&& \hspace{-0.2cm} \bmat{C_i + \sum_{j=1}^\ell [\Lambda_0]_{ij}W_j & (\rho_i/2) e_1 \\[0.8em] 
				(\rho_i/2) e_1^T & -\gamma_i} \succeq 0, \\[0.5em]
	&&& [\Lambda_0]_{ij} \le 0, \quad  \ j \in [\ell].
\end{alignedat} 
\end{equation}
Since $v_0$ is a feasible solution to \eqref{prog:SDR_OPF_polar1}, it holds that
\[
	v_0^\herm A_iv_0 + b_i^T x_0 \le \pi_i = \tau_i.
\]
Therefore,
there exists $\rho_i,\gamma_i, [\Lambda_0]_{i1},\dots,[\Lambda_0]_{i\ell}  \in \R$, such that $[\Lambda_0]_{ij} \le 0,$ for all $j \in [\ell]$, $v_0^\herm  A_iv_0+ b_i^T x_0\le \rho_i+\gamma_i,$ and 
\begin{align}	\label{eq:help_proof_feas1}
		\bmat{C_i + \sum_{j=1}^\ell [\Lambda_0]_{ij}W_j & (\rho_i/2) e_1 \\[0.8em]  
				(\rho_i/2) e_1^T & -\gamma_i} \succeq 0.
\end{align}
\noindent Since  $v_0^\herm A_iv_0 + b_i^T x_0\le \rho_i+\gamma_i,$ it follows readily that  
\begin{align}	\label{eq:help_proof_feas2}
		(v_0^\herm A_iv_0 + b_i^T x_0)e_1e_1^T  \ \preceq \ (\rho_i+\gamma_i) e_1e_1^T.
\end{align}	

\noindent Subtracting $C_i + \sum_{j=1}^\ell [\Lambda_0]_{ij}W_j$ from both sides of \eqref{eq:help_proof_feas2}, we observe that its left-hand side becomes equal to \eqref{eq:help_proof_feas}. Therefore, it suffices to show that 
\begin{align}	\label{eq:help_proof_feas3}
		 (\rho_i+\gamma_i) e_1e_1^T-C_i - \sum\limits_{j=1}^\ell [\Lambda_0]_{ij}W_j \ \preceq \ 0,
\end{align}
as this implies that \eqref{eq:help_proof_feas} holds. Using Lemma 3.1 in \cite{kim2003second},  
the positive semidefiniteness constraint \eqref{eq:help_proof_feas1} can be equivalently described by the following set of conditions:
\begin{equation}\label{eq:cond}
\begin{aligned} 
\begin{cases}
	&C_i +\displaystyle\sum_{j=1}^\ell [\Lambda_0]_{ij}W_j \succeq 0,  \\
	&\gamma_i \le 0, \\
	&\gamma_i\Big(C_i + \displaystyle\sum_{j=1}^\ell [\Lambda_0]_{ij}W_j\Big) + (\rho_i^2/4) e_1e_1^T \preceq 0. 
\end{cases}
\end{aligned}
\end{equation}

\

\noindent Using the equivalent conditions in \eqref{eq:cond}, we verify the satisfaction of the desired  inequality in \eqref{eq:help_proof_feas3} in two separate cases: $\gamma_i = 0$ and $\gamma_i < 0$. 

\

\noindent \emph{Case 1:} \ Let  $\gamma_i = 0$. Then, \eqref{eq:help_proof_feas1} implies that $\rho_i = 0$. This follows from the fact that  for any positive semidefinite matrix $X$,  if the diagonal entry satisfies  $[X]_{mm}=0$, then all off-diagonal entries  in the correspoding row and column also  satisfy $[X]_{ml} =[X]_{lm}= 0,$ for all $l \in [n]$.
In this case, \eqref{eq:help_proof_feas3} becomes equivalent to the first condition in \eqref{eq:cond}, and we are done.

\

\noindent \emph{Case 2:} \ Let  $\gamma_i < 0$. By rearranging terms, the third condition in \eqref{eq:cond} holds if and only if  
\[
	-C_i-\sum_{j=1}^\ell [\Lambda_0]_{ij} W_j \ \preceq \ \frac{\rho_i^2}{4\gamma_i} e_1e_1^T.
\]
It follows that \eqref{eq:help_proof_feas3} holds if 
\begin{align}\label{eq:to_show}
	\rho_i+\gamma_i+ \frac{\rho_i^2}{4\gamma_i}  \le 0.
\end{align}
For a matrix $X \in \Hn$, let $\eta(X) \in \R^{n}$ be the vector of eigenvalues of $X$ arranged in nonincreasing order. According to 
Corollary 7.7.4 in \cite{horn2012matrix}, it holds that $\eta(X) \ge \eta(Y)$  for any pair of matrices $X,Y \in \Hn$ that satisfy $X \succeq Y$. It follows that the third condition in \eqref{eq:cond} implies that
\[
	 \rho_i^2 \le -4\gamma_i \eta_{\max},
\] 
where $\eta_{\max}$ is the largest eigenvalue of $C_i+\sum_{j=1}^\ell [\Lambda_0]_{ij} W_j.$ As this is a positive semidefinite matrix, we must also have that $\eta_{\max}$ is nonnegative.
Therefore, 
$\rho_i \le 2\sqrt{-\gamma_i\eta_{\max}}$.
Using this inequality, we can bound the left hand-side of \eqref{eq:to_show} from above. And therefore, it suffices to show that
\begin{align}\label{eq:to_show_final}
	2\sqrt{-\gamma_i\eta_{\max}}+\gamma_i - \eta_{\max} \ \le \ 0,
\end{align}
as this implies that \eqref{eq:to_show} holds. Indeed, consider the function $f:\R_+\times \R_{++} \rightarrow \R$ given by $f(x,y) = -y+2\sqrt{yx} - x$, where $\R_+$ and $\R_{++}$ denote the sets of nonnegative and positive real numbers, respectively. To complete the proof, notice that the function $f$ has a maximum value of zero over its domain.  \qed

\subsection{Proof of Proposition \ref{prop:Algorithm_properties}} \label{app:proof_iterative}
\noindent For notational brevity, let us first define the matrices
\[
	F_i(Z,x,V,\Lambda) \coloneqq H_i(V,Z)- C_i +(b_i^Tx)e_1e_1^T - \sum_{j=1}^\ell [\Lambda]_{ij} W_j,
\]
for each $i \in [m]$. In addition, define
  $$J(V,V_t)  \coloneqq  \tr(MH_0(V,V_t)).$$ 
The proof is by induction on $t$.

\

\noindent\textit{Base step:} Let $t=0$. By Propostion \ref{prop:feasibility}, there exists $x_0 \in \R^{2n}$ and $\Lambda_0 \in \R^{m \times \ell}$ such that $(x_0,V_0,\Lambda_0) \in \Feas{V_0}$. Hence, $\Feas{V_0} \neq \emptyset$, and we established property \eqref{prop:res_1} for the base step. Let 
\begin{align*}
(x_1,V_1,\Lambda_1) \   \in    \underset{(x,V,\Lambda) \in \Feas{V_0}}{\argmin} \   \tr(MH_0(V,V_0)).
\end{align*}
It follows that 
\[
	J(V_1,V_0)   \le  J(V_0,V_0)  =  \tr(MV_0^\herm A_0V_0),
\]
since $(x_0,V_0,\Lambda_0) \in \Feas{V_0}$ is also feasible.
It remains to show that $\tr(MV_1^\herm A_0V_1) \le J(V_1,V_0)$ in order to establish \eqref{prop:res_2} for the base step.
Recall property \eqref{lemma:properties_Hi_b} of Lemma \ref{lemma:properties_Hi}. Letting $V =V_1$ and $Z = V_0$, we  obtain the desired inequality.

\

\noindent \textit{Induction step:} 
Let $t=s$, and suppose  that $\Feas{V_j} \neq \emptyset$ for all $j < s$. We must show that $\Feas{V_s}\neq \emptyset$ in order to establish \eqref{prop:res_1}. To do so, we use the assumption that $\Feas{V_{s-1}} \neq \emptyset$, together with Lemma \ref{lemma:properties_Hi}\eqref{lemma:properties_Hi_a}  to show that $(x_s,V_s,\Lambda_s)  \in \Feas{V_s}$.   Let $(x_{s}, V_{s}, \Lambda_{s}) \in \Fcal(V_{s-1})$ be an optimal solution to \eqref{eq:iter_alg}. This solution is guaranteed to exist since $\Feas{V_{s-1}}\neq \emptyset$ by assumption.  Thus, $Ex_s \le f$, $\Lambda_s\le 0$, and 
\begin{equation}
\begin{aligned}\label{eq:help_prop2_1}
 F_i(V_{s-1},x_s,V_s,\Lambda_s) \preceq 0, \quad \forall  \ i \in [m].
 \end{aligned}
\end{equation}
Fix an $i \in [m].$ By using Lemma \ref{lemma:properties_Hi}\eqref{lemma:properties_Hi_a} at $V = V_{s}$ and $Z = V_{s-1}$, we obtain  
\begin{align}\label{eq:help_prop2_2}
 	V_{s}^\herm A_iV_{s} \ \preceq \ H_i(V_s,V_{s-1}).
\end{align}
Adding $- C_i +(b_i^T x_s)e_1e_1^T - \sum_{j=1}^\ell [\Lambda_s]_{ij} W_j $ 
to both sides of \eqref{eq:help_prop2_2}, yields 
\begin{alignat}{4}\label{eq:help_prop2_3}
V_{s}^\herm A_iV_{s} - C_i +(b_i^T x_s)e_1e_1^T - &\sum_{j=1}^\ell [\Lambda_s]_{ij} W_j 
	 \\ &\preceq F_i(V_{s-1},x_{s},V_s,\Lambda_s)  \preceq 0. \nonumber
\end{alignat}
where the last relation follows from \eqref{eq:help_prop2_1}.  Consider now the feasible set $\Feas{V_s}$. We claim that $(x_{s},V_{s},\Lambda_{s}) \in \Feas{V_s}.$ Indeed $Ex_s\le f, \ \Lambda_s \le 0$ and $F_i(V_s,x_s,V_s,\Lambda_s)  \preceq 0$ for all $i \in [m]$, 
where the last relation follows from \eqref{eq:help_prop2_3} and the fact that $H_i(V_s,V_s) = V_{s}^\herm A_iV_{s}$. This completes the induction step for part \eqref{prop:res_1}.

We will now establish the induction step for part  \eqref{prop:res_2}. We must show that $\tr(MV_s^\herm A_0V_s)  \le   \tr(MV_{s-1}^\herm A_0V_{s-1})$ for any  optimal solution $(x_s,V_s,\Lambda_s)$ of \eqref{eq:iter_alg}. Let $(x_{s},V_{s},\Lambda_{s})$ be one such solution.
Since it is optimal, we must have
\[
	J(V_s,V_{s-1}) \le J(V,V_{s-1}), \quad \forall \ (x,V,\Lambda) \in \Feas{V_{s-1}}.
\]
In the step of induction of part \eqref{prop:res_1}, we have shown that $(x_{s-1},V_{s-1},\Lambda_{s-1}) \in \Feas{V_{s-1}}$. Therefore, we obtain in particular
\[
	J(V_s,V_{s-1}) \le J(V_{s-1},V_{s-1}) = \tr(MV_{s-1}^\herm A_0V_{s-1}). 
\]
It remains to show that $\tr(MV_{s}^\herm A_0V_{s}) \le J(V_s,V_{s-1}).$ This follows  by setting $V =V_s$ and $Z = V_{s-1}$ in Lemma \ref{lemma:properties_Hi}\eqref{lemma:properties_Hi_b}. \qed

\section{LMI Reformulation of \ref{prog:two_stage_Linearization}} \label{sec:appendix_Prop1}
For $i =0,1,\dots,m$,  let $B_i = (A_i^+)^{1/2}$ and $N = M^{1/2}$ be the square roots of $A_i^+$ and $M$, respectively.  These matrices are guaranteed to exist since both $A_i^+$ and $M$ are positive semidefinite.
In our reformulation, it will be convenient to write
\begin{align*}	
	\tr(MV^\herm A_0^+V) &= \tr\left((B_0VN)^\herm (B_0VN)\right)\\
	 &= \vec(B_0VN)^\herm \vec(B_0VN),
\end{align*}
 where $\vec(\cdot)$ denotes the linear operator vectorizing matrices by stacking their columns.
Let
\begin{align*}
	L_i(V,V_0) &= H_i(V,V_0)-V^\herm A_i^+V\\
		     &= V_0^\herm A_i^-V  +  V^\herm A_i^-V_0  -  V_0^\herm A_i^-V_0 
\end{align*}
denote the part of $H_i(V,V_0)$ that depends affinely on $V$. Applying the Schur complement formula to the  matrix inequalities in \ref{prog:two_stage_Linearization}, we obtain the following equivalent reformulation of \ref{prog:two_stage_Linearization} as a semidefinite program:

\begin{alignat}{8} \label{prog:two_stage_Linearization1}
	& \underset{}{\text{minimize}} \quad  t+\tr(ML_0(V,V_0))  \nonumber \\
	& \text{subject to} \ \ x \in \R^{2n}, \ V \in \C^{n \times k}, \ \Lambda \in \R^{m\times\ell},\ t\in \R \nonumber  \\[0.5em]
	&  \hspace{-0.15cm} \bmat{-I_n & B_iV \\ V^\herm B_i^\herm & L_i(V,V_0) - C_i +(b_i^Tx)e_1e_1^T - \displaystyle\sum\limits_{j=1}^\ell [\Lambda]_{ij} W_j  }  \preceq 0,  \nonumber \\
	&   \forall \ i \in [m]. \nonumber \\[0.5em]
	&  \bmat{-tI_{nk} & \vec(B_0VN) \\ \vec(B_0V N)^\herm & -1 } \preceq 0,
\nonumber \\[0.5em]
	&      Ex   \le f, \nonumber \\[0.3em]
	&   \Lambda   \le 0. \nonumber 
\end{alignat}

\vspace{-0.5cm}

\section{RAC-OPF with Load Shedding} \label{app:Load_Shedding}

In what follows,  we provide an alternative formulation of the RAC-OPF problem, which allows 
for \emph{load shedding} at non-generator buses in the power network. Here,  a reduction in load relative to the nominal demand level  is penalized according to a suitably chosen \emph{value of lost load} (VOLL).\footnote{For example, the  current VOLL  used within the Midcontinent Independent System Operator (MISO) markets is \$3,500/MWh \cite{miso2017report}.}

In order to develop this generalization of the RAC-OPF problem, we first require some additional notation. Let $\Vg \subseteq \Vcal$ denote the subset of buses connected to generators, and define $\VL \coloneqq \Vcal\setminus\Vg$  as the subset of non-generator buses. Also, let $\nG:=|\Vg|$ and $\nL:=|\VL|$ denote the number of generator and non-generator buses, respectively. 
Define  the complex load that is \emph{shed} at each non-generator bus $i \in \VL$ by $\lambda_i(\bm{\xi})$,  where $\lambda_i \in \Lb^2_{k,1}$ is a recourse function determined by the ISO for each load $i \in  \VL$. The power balance equation at each non-generator bus $i \in \VL$ can therefore be expressed as
\[
\lambda_i(\bm{\xi}) -  v(\bm{\xi})^\herm S_i v(\bm{\xi}) = \dnom_i.
\]
Finally, letting $\beta \in \R_+$ denote the VOLL, we arrive at the following reformulation of the RAC-OPF problem with load shedding:
\begin{alignat}{4}
		& \underset{}{\text{minimize}}   
		&& & &  \hspace{-1.7in} \Exp \left[ \sum_{i\in \Vg} \alpha_i \real\{g_i(\bm{\xi})\}  +  \beta \sum_{i \in \VL} \real\{\lambda_i(\bm{\xi})\} + \imag\{\lambda_i(\bm{\xi})\}  \right] \nonumber  \\[0.5em]
		& \text{subject to}  \ \    
		&& & &   \qquad \qquad \qquad \quad \ \  \mathllap{ g^0 \in \C^{\nG}, \ \  g \in \Lb^2_{k,\nG}, \ \ v \in \Lb^2_{k,n}, \ \ \lambda \in\Lb^2_{k,\nL} }     \nonumber \\[0.5em]
		&&& g_i^{\min}  \le  g_i^0 \le  g_i^{\max}, &&    i\in  \Vg    \nonumber \\[0.5em]
		&&& \underline{g}_i(\xi)   \le  g_i(\xi)   \le  \overline{g}_i(\xi),   &&        i\in  \Vg 
				\makebox[48pt][r] {\smash{\raisebox{-4.1\baselineskip}{$\left.\rule{0pt}{4.65\baselineskip}\right\}\forall \ \xi \in \Xi.$}}}
		 \nonumber \\[0.5em]
		&&&  r_i^{\min} \le  g_i(\xi) - g^0_i  \le  r_i^{\max}, &&   i\in  \Vg  \nonumber \\[0.5em]
		&&&  g_i(\xi) -v(\xi)^HS_iv(\xi)  =  \dnom_i,  \quad &&   i\in  \Vg   \hspace{.2in}   \nonumber \\[0.5em]
		&&&  v_i^{\min}  \le  |v_i(\xi)|  \le  v_i^{\max},  && i\in  \Vg  \nonumber \\[0.5em]
		&&&  |v(\xi)^H P_{ij} v(\xi)|  \le   \ell_{ij}^{\max}, && \hspace{-0.63cm}  (i,j)\in \Ecal    \nonumber  \\[0.5em]
			&&&    0 \leq \lambda_i(\xi) \le \dnom_i, &&  i  \in \VL  \hspace{.2in}   \nonumber \\[0.5em]
			&&&  \lambda_i(\xi) -v(\xi)^H S_iv(\xi)  =  \dnom_i,  \quad &&   i\in  \VL  \hspace{.2in}   \nonumber
\end{alignat}

\section{Table of Widely-Used Variables and Symbols}\label{app:table_notation}
\begin{table}[H]
\blueN{\setlength{\tabcolsep}{2.50pt}
\renewcommand{\arraystretch}{1.12} 
\begin{tabular*}{3.65in}{lllcl} 
\toprule 
\multicolumn{1}{c}{Var/Symb} & & Space & & Description  \\
\cmidrule{1-1} \cmidrule{3-3} \cmidrule{5-5} 
$n$ & & $\Nb$  &  & number of network buses   \\    
$m$ & & $\Nb$  &  & number of constraints of RAC-OPF   \\
$k$ & & $\Nb$  &  & ambient dimension of uncertainty set $\Xi$   \\
$\ell$ & & $\Nb$  &  & number of ellipsoids describing  $\Xi$   \\    
$\Gcal = (\Vcal,\Ecal)$ & &  \ --  &  & graph describing the power network   \\      
$\Vcal$ & &  \ --  &  & set of vertices indexing the transmission buses   \\      
$\Ecal$ & &  \ --  &  & set of edges indexing the transmission lines   \\ 
$Y$ & &  $\C^{n\times n}$  &  & bus admittance matrix   \\ 
$\bm{\xi}$ & &   \ --   &  & random vector     \\ 
$\xi$ & &   $\R^{k}$   &  & a realization of the random vector $\bm{\xi}$    \\
$\Xi$ & &  \ --     &  & support of the random vector $\bm{\xi}$    \\
$\mu$ & &   $\R^{k}$   &  & first-order moment of $\bm{\xi}$    \\  
$M$ & &   $\R^{k\times k}$    &  & second-order moment of $\bm{\xi}$   \\     
$\ell_{ij}^{\max}$ & & $\R$  &  & active power flow capacity of transmission line $(i,j)$  \\  
$\alpha_i$ & & $\R_+$ & &  marginal cost of active power generation at bus $i$    \\     
$\beta$ & & $\R_+$ & &  value of lost load (VOLL)  \\                              
$v$ & &  $\Lcal^{2}_{k,n}$  &  & real-time bus voltage phasor   \\  
$g^0_i$  & & $\C$ & & day-ahead dispatch of generator $i$  \\               
$g_i$ & & $\Lb^2_{k,1}$ & &  real-time dispatch of generator $i$     \\
$g_i^{\max}$ & & $\C$ & & nominal maximum capacity of generator $i$   \\
$g_i^{\min}$& & $\C$ & &  nominal minimum capacity of generator $i$   \\
$\overline{g}_i$ & & $\Lb^2_{k,1}$ & & real-time maximum capacity of generator $i$   \\
$\underline{g}_i$& & $\Lb^2_{k,1}$ & &  real-time minimum capacity of generator $i$   \\
$\dnom_i$  &  & $\C$  & &  real-time demand at bus $i$ \\ 
$r_i^{\min}$& & $\C$ & &  ramp-down limit of generator $i$  \\
$r_i^{\max}$& & $\C$ & &  ramp-up limit of generator $i$   \\
$\gamma_i^{\max}$&  & $\C$  & &  guaranteed maximum available power supply at bus $i$\\ 
$\gamma_i^{\min}$&  & $\C$  & &  guaranteed minimum available power supply at bus $i$\\  \bottomrule
\end{tabular*}
\vspace{-0.1cm}}
\end{table}

\end{appendices}

\end{document}